\documentclass[10pt, a4paper]{amsart}
\usepackage[dvips]{epsfig}
\usepackage{amsmath}
\usepackage{amssymb}
\usepackage{amsfonts}
\usepackage{dsfont}
\usepackage{bbm}
\usepackage{amsthm}
\usepackage{amsbsy}
\usepackage{amsgen}
\usepackage{amscd}
\usepackage{amsopn}
\usepackage{amstext}
\usepackage{amsxtra}
\usepackage[utf8]{inputenc}
\usepackage{enumerate}
\usepackage{hyperref}
\usepackage{lineno}
\usepackage{marginnote}
\usepackage{verbatim}
\usepackage{color}
\usepackage{calrsfs}
\usepackage{times, graphicx}
\usepackage{eso-pic}
\usepackage{lastpage}
\DeclareMathAlphabet{\pazocal}{OMS}{zplm}{m}{n}
\usepackage[foot]{amsaddr}
\usepackage{geometry}

\newtheorem{theorem}{Theorem}[section]
\newtheorem{lemma}[theorem]{Lemma}

\theoremstyle{definition}

\newtheorem{proposition}[theorem]{Proposition}
\newtheorem{corollary}[theorem]{Corollary}

\theoremstyle{remark}

\numberwithin{equation}{section}

\newcommand{\cal}{\mathcal}
\allowdisplaybreaks

\newcommand {\E} {\mathbb{E}}

\newcommand {\R} {\mathbb{R}}

\newcommand{\id}{\mathds{1}}

\begin{document}

\title[On the correlation between critical points and critical values]{On the correlation between critical points and critical values for random spherical harmonics}

\author{Valentina Cammarota}
\address{Department of Statistics, Sapienza University of Rome.}

\email{valentina.cammarota@uniroma1.it}
\thanks{The author V. C. has received funding from Sapienza University research project RM120172B80031BE, {\it Geometry of Random Fields}}

\author{Anna Paola Todino}
\address{Department of Mathematical Sciences, Politecnico di Torino.}
\email{anna.todino@polito.it}
\thanks{The author A.P. T. was partially supported by Progetto di Eccellenza, Dipartimento di Scienze Matematiche, Politecnico di Torino, CUP: E11G18000350001 and by GNAMPA-INdAM (\textit{project: Stime asintotiche: principi di invarianza e grandi deviazioni})}





\begin{abstract}
We study the correlation between the total number of critical points of random spherical harmonics and the number of critical points with value in any interval $I \subset \mathbb{R}$. We show that the correlation is asymptotically zero, while the partial correlation, after controlling the random $L^2$-norm on the sphere of the eigenfunctions, is asymptotically one. Our findings complement the results obtained by Wigman (2012) and Marinucci and Rossi (2021) on the correlation between nodal and boundary length of random spherical harmonics.  
\end{abstract}

\maketitle


\section{Introduction and Main Result}

\subsection{Random spherical harmonics}

Let $\mathbb{S}^2$ be the unit $2$-dimensional sphere and $\Delta_{\mathbb{S}^2}$ be the Laplace-Beltrami operator on $\mathbb{S}^2$. The spectrum of $\Delta_{\mathbb{S}^2}$ consists of the numbers $\lambda_\ell =\ell(\ell+1)$ with $\ell =1, 2, \dots$, and the eigenspace corresponding to $\lambda_\ell$ is the $(2 \ell+1)$-dimensional linear space of degree $\ell$ spherical harmonics. For $\ell >0$ let $\{Y_{\ell m}(\cdot)\}_{m=-\ell,\dots,\ell}$ be an arbitrary $L^2$-orthonormal basis of real valued spherical harmonics satisfying 
$$\Delta_{\mathbb{S}^2} Y_{\ell m} +\lambda_\ell Y_{\ell m}=0.$$
On $\mathbb{S}^2$ we consider a family of Gaussian random fields, defined of a suitable probability space $(\Omega, \cal{F}, \mathbb{P})$,   
 \begin{equation} \label{ttt}f_\ell(x)=\frac{\sqrt{4 \pi }}{\sqrt{2\ell+1}} \sum_{m=-\ell}^\ell a_{\ell m} Y_{\ell m}(x), \end{equation}
 where the coefficients $\{a_{\ell m}\}_{m=-\ell,\dots,\ell}$  are independent standard Gaussian with zero mean and unit variance. 
 The standardization in \eqref{ttt} is such that ${\rm Var}(f_{\ell}(x))=1$, and the law of the process $\{f_{\ell}(\cdot)\}$ is invariant with respect to the choice of the $L^2$-orthonormal basis $\{Y_{\ell m}\}$. The random fields $\{f_\ell(x):  x\in \mathbb{S}^2 \}$ are isotropic centred Gaussian with covariance function given by
$$\mathbb{E}[f_\ell(x)f_\ell(y)]=P_\ell(\cos d(x,y)),$$
denoting with $P_\ell$ the Legendre polynomial and $d(x,y)=\arccos \langle x,y\rangle$ the geodesic distance on the sphere.

In this paper, we focus on the critical points and critical values of $f_\ell$. Let $I\subseteq \mathbb{R}$ be any interval in the real line and $\nabla$ the covariant gradient on the sphere, the number of critical points of $f_\ell$ with value in $I$ is denoted by
\begin{equation*}
\mathcal{N}^c_{\ell}(I)= \#\{ x\in \mathbb{S}^2: \nabla f_\ell(x)=0, f_\ell(x) \in I \};
\end{equation*}
we denote $\mathcal{N}^c_\ell(u)=\mathcal{N}^c_{\ell}(-\infty,u)$ and $\mathcal{N}^c_\ell=\mathcal{N}^c_{\ell}(\mathbb{R})$ the total number of critical points. In this paper, in particular, we investigate how much the number of critical points characterizes the geometry of the random spherical eigenfunctions, i.e. the behaviour of the excursion sets
$$A_{u}(f_{\ell})=\{x \in \mathbb{S}^2: f_{\ell}(x) \ge u\},$$
for arbitrary levels $u \in \mathbb{R}$.

A number of issues on the geometry of random spherical harmonics has been recently analysed: nodal domains  \cite{nazarov, Log Upper}, length of nodal lines \cite{Wig, MRW}, the excursion area and the
defect \cite{DI, MW, MR2015}, Euler-Poincar\'{e}
characteristic of the excursion sets \cite{CMW-EPC, CM2018}, mass
equidistribution \cite{han}, critical radius \cite{feng}. These and other
geometric features have also been intensively studied for random eigenfunctions on other manifolds such
as the torus (Arithmetic Random Waves) and the plane (Berry's Random Waves
model), see e.g.  \cite{Berry 1977, BMW, KKW, MPRW2015, Buckley, BCW, estradeleon, Cammarota2018, CMR, PV, Vid}; \cite{BMW, Todino1, Todino2} for fluctuations over subdomains of the torus and of the sphere, \cite{granville} for the analysis of mass equidistributions; and \cite{Rudnick, RudnickYesha} for nodal
intersections, to list only some of the recent contributions.

\subsection{Critical values}
In \cite{CMW14} it has been shown that, for every interval $I \subseteq \mathbb{R}$, as $\ell \to \infty$, the expected number of critical points with value in $I$ behaves like 
$$\mathbb{E}[\mathcal{N}^c_\ell(I)]=\frac{2}{\sqrt{3}} \ell^2\int_I \frac{\sqrt{3}}{\sqrt{8\pi}}(2e^{-t^2}+t^2-1)e^{-\frac{t^2}{2}} dt +O(1),$$
where here (and later) the constant in the $O(\cdot)$ term is universal, i.e. the integral of the error term on any interval $I$ is uniformly bounded by its value when $I=\mathbb{R}$.  
The investigation of the asymptotic variance is more challenging and in \cite[Theorem 1.2]{CMW14} it has been shown that
\begin{equation} \label{varNu}
\text{Var}(\mathcal{N}^c_\ell(I))=\ell^3 [\nu^c(I)]^2+O(\ell^{5/2}),
\end{equation}
where 
$$\nu^c(I)=   \int_I \frac{1}{\sqrt{8\pi}} [2-6 t^2 - e^{t^2} (1-4t^2 + t^4 ) ]e^{-\frac{3 }{2} t^2}\,dt.$$
Similar results hold for the number of extrema and saddles. 
\subsection{Critical points}
When considering the total number of critical points, i.e. $I=\mathbb{R}$, we immediately obtain that 
$$\mathbb{E}[\mathcal{N}^c_\ell]=\frac{2}{\sqrt{3}} \ell^2 +O(1),$$
whereas the leading term in \eqref{varNu} vanishes and \cite[Theorem 1.1]{CW15} establishes that as $\ell \to \infty$,
\begin{equation}\label{varN}
\text{Var}(\mathcal{N}_\ell^c)= \frac{1}{3^3 \pi^2}\ell^2 \log \ell+O(\ell^2).
\end{equation}
\subsection{Interpretation in terms of Wiener chaoses} 
These results can be interpreted in terms of the $L^2(\Omega)$ expansion of critical points into Wiener chaoses, see e.g. \cite{CM2018b}, which are orthogonal spaces spanned by Hermite polynomials. First of all, we recall that the Hermite polynomials $H_q(x)$ are defined by $H_0(x)=1$, and for $q=2,3,\dots$ 
\begin{equation*}
H_q(x)=(-1)^q \frac{1}{\phi(x)} \frac{d^q\phi(x)}{dx^q},
\end{equation*}
with $ \phi(x)=\frac{1}{\sqrt{2\pi}}e^{-x^2/2}$. We consider the Wiener chaos expansion 
\begin{align}\label{Proj}
\mathcal{N}^c_\ell(I)=\sum_{q=0}^{\infty} \mathcal{N}^c_\ell(I)[q],
\end{align}
where $\mathcal{N}^c_\ell(I)[q]$ denotes the projection of $\mathcal{N}^c_\ell(I)$ on the $q$-order chaos component that is the space generated by the $L^2$-completion of linear combinations of the form 
$$H_{q_1}(\xi_1) \cdot H_{q_2}(\xi_2) \cdots H_{q_k}(\xi_k), \hspace{1cm} k \ge 1,$$ 
with $q_i \in \mathbb{N}$ such that $q_1+\cdots+q_k=q$, and $(\xi_1, \dots, \xi_k)$ standard real Gaussian vector. 

 It results that (after centring) a single term dominates the $L^2(\Omega)$ expansion in (\ref{Proj}).  We define the random variables 
\begin{equation*}
h_{\ell,q}:= \int_{\mathbb{S}^2} H_q(f_\ell(x))\,dx
\end{equation*}
called sample polyspectra, see i.e.  \cite{MW, MW2014, PTRF09}. We have that 
$$\text{Var}(h_{\ell,2})=(4\pi)^2 \frac{2}{2\ell+1}, \hspace{1.5cm} \text{Var}(h_{\ell,4})=576 \frac{\log \ell}{\ell^2}+O(\ell^{-2}),$$
and, for $q=3$ and $q\geq 5$, 
$$ \text{Var}(h_{\ell,q})= \frac{c_q}{\ell^2}+o(\ell^{-2}), \hspace{1.5cm} c_q:=\int_0^{\infty} J_0(\psi)^q\psi\,d\psi,$$ and $J_0(\cdot)$ is the Bessel function of order zero. Note that the coefficient $c_3$ can be calculated in the closed form $c_3=\frac{2}{\pi \sqrt{3}}$ (see eq. (2.12.42.15) in \cite{PBM}). 

These results, and equations (\ref{varNu}) and (\ref{varN}), suggest that the asymptotic behaviour of the total number of critical points is dominated by the projection into the fourth chaotic component, which can be expressed by the integral of $h_{\ell,4}$;  whereas the number of critical values in $I$ is dominated by the projection into the second chaotic component, which can be expressed by $h_{\ell,2}$. Indeed let us introduce the random variables  
$$\mathcal{S}_\ell(I)=\frac{\lambda_\ell}{2}  \nu^c(I) \frac{1}{2\pi} \int_{\mathbb{S}^2} H_2(f_\ell(x))dx,$$
and
$$\mathcal{F}_\ell=- \frac{\lambda_\ell}{2^33^2\sqrt{3} \pi}  \int_{\mathbb{S}^2} H_4(f_\ell(x))dx.$$
In \cite{CM2018b} it has been established that, as $\ell \to \infty$, and for $I \subset \mathbb{R}$ and such that $\nu^c(I) \ne 0$, 
$${\mathcal{N}}_{\ell}^{c}(I) -\mathbb{E}[{\mathcal{N}}_{\ell}^{c}(I)] ={\mathcal{N}}_{\ell}^{c}(I)[2] + R_{\ell}(I),$$
with $\mathbb{E}[R^2_{\ell}(I)]=o(\ell^3)$ uniformly over $I$, and where ${\mathcal{N}}_{\ell}^{c}(I)[2]=\mathcal{S}_\ell(I)$. As a consequence the total number of critical values is fully correlated in the limit with $\mathcal{S}_\ell(I)$, i.e. as $\ell \to \infty$     
\begin{equation}\label{correlation1}
\text{Corr}(\mathcal{N}_\ell^c(I), \mathcal{S}_\ell(I))=\frac{\text{Cov}(\mathcal{N}_\ell^c(I), \mathcal{S}_\ell(I))}{\sqrt{ \text{Var}(\mathcal{N}_\ell^c(I)) \text{Var}(\mathcal{S}_\ell(I))}} \to 1.  
\end{equation}
Subsequently in \cite{CM2019} it has been shown that 
\begin{equation*}
	\mathcal{N}_{\ell }^{c}-\mathbb{E}\left[ \mathcal{N}_{\ell }^{c}\right]
	={\mathcal{N}}_{\ell}^{c}[4]+o_{\mathbb{P}}(\sqrt{\ell^2 \log \ell}),
	\end{equation*} 
where ${\mathcal{N}}_{\ell}^{c}[4]=\mathcal{F}_\ell$, and, in general, for e sequence of random variables $X_\ell$ and a sequence of real numbers $a_\ell$, the notation  $X_\ell=o_{\mathbb{P}}(a_{\ell})$ means that $X_\ell/ a_{\ell}$ converges to zero in probability as $\ell \to \infty$. Hence the total number of critical points is fully correlated in the limit with $\mathcal{F}_\ell$
\begin{equation*}
\text{Corr}(\mathcal{N}_\ell^c, \mathcal{F}_\ell)=\frac{\text{Cov}(\mathcal{N}_\ell^c, \mathcal{F}_\ell)}{\sqrt{\text{Var}(\mathcal{N}_\ell^c)\text{Var}(\mathcal{F}_\ell)}} \to 1. 
\end{equation*}
An important consequence of the results in \cite{CM2018b, CM2019}, is that, while the computation of the number of critical points and critical values via Kac-Rice formula (see \cite{adlertaylor} and Section \ref{KRF} below) requires the evaluation of gradient and Hessian fields, the dominant term of $\mathcal{N}_{\ell }^{c}$ and $\mathcal{N}_{\ell }^{c}(I)$ depends, in the high frequency limit, only on the second-order and fourth-order Hermite polynomials evaluated at the eigenfunctions $f_{\ell}$, i.e. only on $h_{\ell,4}$ and $h_{\ell,2}$ respectively. 
Moreover 
\begin{align*}
 h_{\ell,2}= \int_{\mathbb{S}^2} f^2_\ell(x) \,dx-4\pi =\frac{4\pi}{2\ell+1} \sum_{m=-\ell}^{\ell} |a_{\ell m}|^2-\mathbb{E}|a_{\ell m}|^2  
\end{align*}
is proportional to a sum of independent and identically distributed random variables with zero mean and finite variance and, as a simple corollary, this implies a quantitative Central Limit Theorem for $\mathcal{N}_{\ell }^{c}(I)$. Similarly for $\mathcal{N}_{\ell }^{c}$: the limiting distribution of $h_{\ell ;4}$ was  studied
in \cite{MW2014}, where it is shown that a quantitative version of the
Central Limit Theorem holds for $h_{\ell ;4}$.

\subsection{Main results} 
The main result in this paper is the characterization of the correlation structure between the critical points and the critical values in any interval $I$. More precisely, we prove that the correlation between $\mathcal{N}^c_\ell(I)$ and $\mathcal{N}_\ell^c$ is asymptotically zero when $I \ne \mathbb{R}$ and $\nu^c(I) \ne 0$, while the partial correlation, after controlling for the random $L^2$-norm on the sphere of the eigenfunctions, is asymptotically one. The proof follows the lines of \cite{MR} where an analogous result is obtained for the correlation between nodal length and boundary length of excursion sets. To do so, we first recall the definition of partial correlation coefficient between two random variables $X_i$, $i=1,2$, with respect to a random variable $Z$
\begin{equation}
\text{Corr}_Z(X_1,X_2)= \text{Corr}(X_1^*,X_2^*),
\end{equation}
where the random variables $X_i^*$ are defined by  
\begin{equation*}\label{X}
X_i^*:= (X_i-\mathbb{E}[X_i])-\frac{\text{Cov}(X_i,Z)}{\text{Var}(Z)}(Z-\mathbb{E}[Z]).
\end{equation*}
In our context the random variables involved are 
$$X_1=\mathcal{N}_\ell^c, \hspace{1cm }X_2=\mathcal{N}_\ell^c(I), \hspace{1cm} Z=||f_\ell(x)||^2_{L^2(\mathbb{S}^2)},$$ and so the partial correlation coefficient measures the linear dependence between $\mathcal{N}_\ell^c$ and $\mathcal{N}_\ell^c(I)$ after getting rid of the components depending on the random $L^2$-norm of the eigenfunctions $f_{\ell}$. Note that 
\begin{align*}
Z-\mathbb{E}[Z]&= ||f_\ell(x)||^2_{L^2(\mathbb{S}^2)}-\mathbb{E}||f_\ell(x)||^2_{L^2(\mathbb{S}^2)}\\
& = \int_{\mathbb{S}^2} f_\ell(x)^2 \,dx-4\pi =\int_{\mathbb{S}^2} H_2(f_\ell(x))\,dx =\frac{4\pi}{2\ell+1} \sum_{m=-\ell}^{\ell} |a_{\ell m}|^2-\mathbb{E}|a_{\ell m}|^2. 
\end{align*}
Assuming the subset $I_1 \subseteq \mathbb{R}$ is such that $\nu^c(I_1)=0$, we prove the following: 
\begin{equation*}
\lim_{\ell \to \infty} \text{Corr}(\mathcal{N}_\ell^c(I_1),\mathcal{N}_\ell^c(I_2))= \begin{cases}0 & \text{if  }  \nu^c(I_2) \ne 0, \\ 1 & \text{if  } \nu^c(I_2) =0, \end{cases} 
\end{equation*}
and, for every $I_1, I_2 \subseteq \mathbb{R}$,
\begin{equation*}
\lim_{\ell \to \infty} \text{Corr}_{||f_\ell(x)||^2_{L^2(\mathbb{S}^2)}}(\mathcal{N}_\ell^c(I_1),\mathcal{N}_\ell^c(I_2))=1.
\end{equation*}
We state our main result taking in particular $I_1 = \mathbb{R}$.  
\begin{theorem}\label{main theorem}
	For subsets $I \subset \mathbb{R}$ such that $ \nu^c(I)\ne 0$, 
	$$\lim_{\ell \to \infty} {\rm{Corr}}(\mathcal{N}_\ell^c, \mathcal{N}_\ell^c(I)) = 0,$$
	and for every $I \subseteq \mathbb{R}$ 
	$$ \lim_{\ell \to \infty} {\rm{Corr}}_{||f_{\ell} ||^2_{L^2(\mathbb{S}^2)}}(\mathcal{N}_\ell^c, \mathcal{N}_\ell^c(I)) = 1. $$
\end{theorem}
As observed in \cite{MR} in the case of nodal and boundary lengths, a corollary of Theorem \ref{main theorem} is that $\mathcal{N}_\ell^c$ and $\mathcal{N}_\ell^c(I)$ are asymptotically independent, but, when the effect of the sample norm of $f_{\ell}$ is properly subtracted, their joint distribution is completely degenerate and so the behaviour of the fluctuations of $\mathcal{N}_\ell^c(I)$ is fully explained by $\mathcal{N}_\ell^c$, in the high energy limit. More precisely, denoting $\widehat{\mathcal{N}}_\ell^c:= \frac{\mathcal{N}_\ell^c}{\sqrt{\mbox{Var}(\mathcal{N}_\ell^c)}}$ and $\widehat{\mathcal{N}}_\ell^c(I):=\frac{\mathcal{N}_\ell^c(I)}{\sqrt{\mbox{Var}(\mathcal{N}_\ell^c(I))}}$, it is possible to prove that as $\ell \to \infty$, for  $I \subset \mathbb{R}$ such that $ \nu^c(I)\ne 0$, 
$$(\widehat{\mathcal{N}}_\ell^c, \widehat{\mathcal{N}}_\ell^c(I)) \stackrel{law}{\to} (Z_1, Z_2), \hspace{1cm}(\widehat{\mathcal{N}}_\ell^{c*}, \widehat{\mathcal{N}}_\ell^{c*}(I)) \stackrel{law}{\to} (Z, Z),$$
where $(Z_1, Z_2)$ is a bivariate vector of standard independent Gaussian variables, and $Z$ denotes a standard Gaussian variable.

\subsection{Discussion and Further Result}\label{EPC0}
In \cite{WigS}, see formula (1.9), Wigman has shown that the length of the level curves is asymptotically fully correlated. Our results fit in the framework of the literature which has investigated the relationship between geometric functionals of excursion sets of random spherical harmonics at different levels as in \cite{WigS, CM2019, MR}. Let us recall the definition of the excursion sets of $f_\ell$ at level $u$
$$A_{u}(f_\ell):=\{ x \in \mathbb{S}^2: f_\ell(x)  \ge u \}.$$
The functionals which describe the geometry of such sets are the so called Lipschitz-Killing Curvatures, which correspond to the area, half of the boundary length and the Euler-Poincar\'e characteristic of $A_{u}(f_\ell)$ and are denoted by $\mathcal{L}_2(u,\ell)$, $\mathcal{L}_1(u,\ell)$, $\mathcal{L}_0(u,\ell)$, respectively.

Previous works, see for instance \cite{CM2018, MR2015, MW2014}, show 
that, when $u \ne 0$ (and $u\neq 1,-1$ for the Euler-Poincar\'e characteristic), the three Lipschitz-Killing curvatures are asymptotically fully correlated to $h_{2;\ell}$ in the high frequency limit, namely 
$$\lim_{\ell \to \infty} \mbox{Corr}(\mathcal{L}_k(u,\ell), h_{2;\ell})=1, \hspace{1cm} k=0,1,2.$$ 
Then, we also immediately have that 
$$\lim_{\ell \to \infty} \mbox{Corr}(\mathcal{L}_k(u_1,\ell),\mathcal{L}_k(u_2,\ell))=1, \hspace{1cm} k=0,1,2$$ for all $u_1,u_2\neq 0$ (and $u\neq 1,-1$ for the Euler-Poincar\'e characteristic). 

Formula \eqref{correlation1} entails that the number of critical values is perfectly correlated, as $\ell \to \infty$, with the area, the Euler-Poincar\'e characteristic and the boundary length at any nonzero levels.
Hence, for $u\ne 0$ (and $u\neq 1,-1$ for the Euler-Poincar\'e characteristic), 
$$\lim_{\ell \to \infty} \mbox{Corr}(\mathcal{L}_k(u,\ell), \mathcal{N}_\ell(u,\infty))=1, \hspace{1cm} k=0,1,2.$$

When the nodal case is considered ($u=0$) the leading term corresponding to $h_{2;\ell}$ of all these geometrical functionals vanishes and the asymptotic behaviour is different. In \cite{MR} the correlation between the nodal length $\mathcal{L}_1(0,\ell)$ and the boundary length $\mathcal{L}_1(u,\ell)$, $u \ne 0$, is investigated; it results that 
$$ \lim_{\ell \to \infty} \mbox{Corr}(\mathcal{L}_1(0,\ell), h_{4;\ell})=1,$$ 
and, for $u \ne 0$,  
$$\lim_{\ell \to \infty} \mbox{Corr}(\mathcal{L}_1(0,\ell), \mathcal{L}_1(u ,\ell))=0,$$ while, after removing the effect of the norm $||f_\ell(x)||_{L^2(\mathbb{S}^2)}$, for any $u \in \mathbb{R}$, it holds that 
$$\lim_{\ell \to \infty}\mbox{Corr}_{||f_\ell(x)||_{L^2(\mathbb{S}^2)}}(\mathcal{L}_1(0,\ell), \mathcal{L}_1(u,\ell))=1.$$ 

Theorem \ref{main theorem} shows that a similar result holds between critical points and critical values: critical values and critical points are asymptotically independent, hence critical points carry no information about the other geometrical functionals at any non-zero levels. This result is due to the fact that the sample norm dominates the behaviour of Lipschitz-Killing curvatures of the excursion sets at non-zero levels, when its effect is adequately removed, the behaviour of $\mathcal{L}_1(u ,\ell)$ at any level is fully explained by the total number of critical points, in the high frequency limit. We have, for $u \ne 0$, 
$$\lim_{\ell \to \infty} \mbox{Corr}(\mathcal{L}_k(u,\ell), \mathcal{N}_\ell^c)=0, \hspace{1cm} k=0,1,2,$$
while 
$$\lim_{\ell \to \infty} \mbox{Corr}(\mathcal{L}_1(0,\ell), \mathcal{N}_\ell^c)=1,$$
and for $u \ne 0$ 
$$\lim_{\ell \to \infty} \mbox{Corr}_{||f_\ell(x)||_{L^2(\mathbb{S}^2)}}(\mathcal{L}_1(u,\ell), \mathcal{N}_\ell^c)=1.$$

A further result of this paper concerns the Euler-Poincar\'e  characteristic. For this geometrical functional a result analogous to Theorem \ref{main theorem} does not hold. In Section \ref{EPCsection} we prove that at level 0 also the fourth chaotic component of the Wiener chaos expansion of the Euler-Poincar\'e  characteristic vanishes. 
In \cite{CMW-EPC, CM2018} it is shown that
$$\mbox{Var}(\mathcal{L}_0(u,\ell))=\frac{\ell^{3}}{4} \left[H_1(u) H_2(u) \frac{e^{-u^2/2}}{\sqrt{2 \pi}} \right]^2+O(\ell^2 \log^2 \ell ),$$
and that the high frequency behaviour of $\mathcal{L}_0(u,\ell)$ is dominated by the projection onto the second order chaos 
 $$\mathcal{L}_0(u,\ell)- \mathbb{E}[\mathcal{L}_0(u,\ell)]= \mathcal{L}_0(u,\ell)[2]+ o_{\mathbb{P}}(\sqrt{\mbox{Var}(\mathcal{L}_0(u,\ell))}),$$
 with 
 $$\mathcal{L}_0(u,\ell)[2]=\frac{\ell^2}{2} \left[H_1(u) H_2(u) \frac{e^{-u^2/2}}{\sqrt{2 \pi}} \right] h_{\ell,2} + R(\ell)$$
 where $\mathbb{E}[R^2(\ell)]=O(\ell^2 \log \ell)$. The projection onto the second order chaos term disappears in the nodal case. However, differently from what happens with nodal length and critical points, the fourth chaotic component is not dominant as well; indeed in Section \ref{EPCsection} we prove that,  in the nodal case
 \begin{proposition}\label{EPC4}
\begin{align*}
\mathcal{L}_0(0,\ell)[4]=0.
\end{align*}
 \end{proposition}

\section{Kac-Rice Formula and $L^{2}$-Convergence} \label{KRF}

By means of Kac-Rice formula, the number of critical points with value in $I$ can be
formally written as
\begin{equation*}
\mathcal{N}_{\ell }^{c}(I)=\int_{\mathbb{S}^{2}}|\text{det}\nabla ^{2}f_{\ell}(x)| \mathbb{I}_{\left\{ f_{\ell }(x)\in I\right\}}
{ \delta} (\nabla f_{\ell }(x))dx,
\end{equation*}
where the identity holds both almost surely (using i.e., the Federer's
coarea formula, see \cite{adlertaylor}), and in the $L^{2}$ sense. 
The validity of this limit in $L^{2}(\Omega )$ was shown in \cite
{CM2018} where it is proved that it is possible to built an approximating sequence of functions 
${\mathcal{N}_{\ell ,\varepsilon }^{c}(I)}$ and establish their convergence both $\omega $-almost surely and in $L^{2}(\Omega )$ to ${\mathcal{N}}%
_{\ell }^{c}(I)$. More precisely, let $\delta _{\varepsilon }:\mathbb{R}%
^{2}\rightarrow \mathbb{R}$ be such that
\begin{equation*}
\delta _{\varepsilon }(z):=\frac{1}{\varepsilon ^{2}}\mathbb{I}_{\{z\in
	\lbrack -\varepsilon /2,\varepsilon /2]^2\}},
\end{equation*}%
and define the approximating sequence
\begin{equation*}
{\mathcal{N}}_{\ell ,\varepsilon }^{c}(I):=\int_{\mathbb{S}^{2}}|\text{det}
\nabla^2f_{\ell }(x)| \mathbb{I}_{\left\{ f_{\ell }(x)\in I\right\}} \delta_{\varepsilon }(\nabla f_{\ell }(x))dx;
\end{equation*}%
it is possible to prove that 
\begin{lemma}
	\label{XXe} For every $\ell \in \mathbb{N}$, we have
	\begin{equation}
	{\mathcal{N}}_{\ell }^{c}(I)=\lim_{\varepsilon \rightarrow 0}{\mathcal{N}}%
	_{\ell ,\varepsilon }^{c}(I),  \label{Xe1}
	\end{equation}%
	where the convergence holds both $\omega $-a.s. and in $L^{2}(\Omega )$.
\end{lemma}

\section{Chaos Expansion}

Following the same approach as given for other geometric functionals in recent papers, 
see e.g. \cite{MPRW2015, CM2018}, we shall start by computing the $L^{2}(\Omega )$ expansion of
critical points into Wiener chaoses, which will lead to 
\begin{equation*}
\mathcal{N}_{\ell }^{c}(I)=\sum_{q=0}^{\infty }\mathcal{N}_{\ell}^{c}(I)[q],  \label{L2exp}
\end{equation*}
where $ \mathcal{N}_{\ell }^{c}(I)[q]$ denotes the chaos-component of order $q$,
or equivalently the projection of $\mathcal{N}_{\ell }^{c}$ on the $q$th order chaos component, which we shall describe below. In order to define and compute more explicitly the chaos components, let
us introduce the standard spherical coordinates $\theta \in [0, \pi]$, $\varphi \in [0, 2 \pi)$ and denote $(\theta_x, \varphi_x)$ the spherical coordinates of $x \in \mathbb{S}^2$, we introduce then the differential operators
\begin{align*}
\partial _{1;x} =\left. \frac{\partial }{\partial \theta }\right\vert
_{\theta =\theta _{x},\varphi =\varphi _{x}} \hspace{1cm} \partial _{2;x}=\left.
\frac{1}{\sin \theta }\frac{\partial }{\partial \varphi }\right\vert
_{\theta =\theta _{x},\varphi =\varphi _{x}}
\end{align*}
\begin{align*}
\partial _{11;x} =\left. \frac{\partial ^{2}}{\partial \theta ^{2}}
\right\vert _{\theta =\theta _{x},\varphi =\varphi _{x}} \hspace{0.4cm} \partial
_{12;x}=\left. \frac{1}{\sin \theta }\frac{\partial^{2}}{\partial \theta
	\partial \varphi }\right\vert_{\theta =\theta _{x}, \varphi =\varphi _{x}} \hspace{0.4cm}
\partial _{22;x}=\left. \frac{1}{\sin ^{2}\theta }\frac{\partial ^{2}
}{\partial \varphi ^{2}}\right\vert _{\theta =\theta _{x},\varphi =\varphi
	_{x}}.
\end{align*}

Recall first that, since $f_{\ell }$ are eigenfunctions of the spherical
Laplacian, for every $x\in \mathbb{S}^{2}$, we can write
\begin{equation}
f_{\ell }(x)=-\frac{\Delta _{\mathbb{S}^{2}}f_{\ell }(x)}{\lambda _{\ell }};  \label{linear_dp}
\end{equation}%
note that at the critical points we have $\Delta_{\mathbb{S}^{2}}f_{\ell
}(x)=\partial _{11}f_{\ell }(x)+\partial _{22}f_{\ell }(x)$, whence their number with
value in $I$ is, by symmetry, given by
\begin{align*}
\mathcal{N}_{\ell }^{c}(I)&=\#\{x\in \mathbb{S}^{2}: \nabla f_{\ell}(x)=0, \; \frac{\partial _{11}f_{\ell }(x)+\partial_{22}f_{\ell}(x)}{\lambda _{\ell }}\in I\}.
\end{align*}
Covariant gradient and Hessian follow the standard definitions, discussed,
for instance, in \cite{CM2018}. Here we simply recall that
\begin{align*}
\nabla f_{\ell }(x) &=(\partial_{1}f_{\ell }(x),\partial_{2}f_{\ell }(x)), \\
\nabla ^{2}f_{\ell }(x) &=\left(
\begin{array}{cc}
\partial _{11}f_{\ell }(x) & \partial _{12}f_{\ell }(x)-\cot \theta
_{x}\partial _{2}f_{\ell }(x) \\
\partial_{12}f_{\ell }(x)-\cot \theta _{x}\partial_{2}f_{\ell }(x) &
\partial_{22}f_{\ell }(x)+\cot \theta_{x}\partial_{1}f_{\ell }(x)
\end{array}
\right), \\
\text{vec}\nabla ^{2}f_{\ell }(x) &=\left(\partial _{11}f_{\ell }(x),\partial
_{12}f_{\ell }(x)-\cot \theta _{x}\partial _{2}f_{\ell }(x),\partial
_{22}f_{\ell }(x)+\cot \theta _{x}\partial _{1}f_{\ell }(x)\right).
\end{align*}
We can then introduce the $5$-dimensional vector $(\nabla f_{\ell }(x),\text{vec}\nabla
^{2}f_{\ell }(x))$; its covariance matrix $\sigma _{\ell }$ is constant with respect to $x$
and it is computed in \cite{CM2018b}. It can be written in the partitioned form
\begin{equation*}
\sigma _{\ell}=\left(
\begin{array}{cc}
a_{\ell } & b_{\ell } \\
b_{\ell }^{T} & c_{\ell }
\end{array}
\right),
\end{equation*}
where the superscript $T$ denotes the conjugate transpose, and
\begin{equation*}
a_{\ell }=\left(
\begin{array}{cc}
\frac{\lambda _{\ell }}{2} & 0 \\
0 & \frac{\lambda _{\ell }}{2}
\end{array}
\right), \hspace{0.5cm}b_{\ell }=\left(
\begin{array}{ccc}
0 & 0 & 0 \\
0 & 0 & 0
\end{array}
\right), \hspace{0.5cm}
c_{\ell }=\frac{\lambda _{\ell }^{2}}{8}\left(
\begin{array}{ccc}
3-\frac{2}{\lambda _{\ell }} & 0 & 1+\frac{2}{\lambda _{\ell }} \\
0 & 1-\frac{2}{\lambda _{\ell }} & 0 \\
1+\frac{2}{\lambda _{\ell }} & 0 & 3-\frac{2}{\lambda _{\ell }}
\end{array}
\right).
\end{equation*}
\noindent Let us recall that the Cholesky decomposition of a
Hermitian positive-definite matrix $A$ takes the form $A=\Lambda \Lambda
^{T},$ where $\Lambda $ is a lower triangular matrix with real and positive
diagonal entries. It is well-known that every Hermitian positive-definite matrix
(and thus also every real-valued symmetric positive-definite matrix) admits
a unique Cholesky decomposition.

By an explicit computation, it is possible to show that the Cholesky
decomposition of $\sigma _{\ell }$ takes the form $\sigma _{\ell }=\Lambda
_{\ell }\Lambda _{\ell }^T$, where
\begin{equation*}
\Lambda _{\ell }=\left(
\begin{array}{ccccc}
\frac{\sqrt{\lambda }_{\ell }}{\sqrt{2}} & 0 & 0 & 0 & 0 \\
0 & \frac{\sqrt{\lambda }_{\ell }}{\sqrt{2}} & 0 & 0 & 0 \\
0 & 0 & \frac{\sqrt{\lambda _{\ell }}\sqrt{3\lambda _{\ell }-2}}{2\sqrt{2}}
& 0 & 0 \\
0 & 0 & 0 & \frac{\sqrt{\lambda _{\ell }}\sqrt{\lambda _{\ell }-2}}{2\sqrt{2}
} & 0 \\
0 & 0 & \frac{\sqrt{\lambda _{\ell }}(\lambda _{\ell }+2)}{2\sqrt{2}\sqrt{
		3\lambda _{\ell }-2}} & 0 & \frac{\lambda _{\ell }\sqrt{{\lambda _{\ell }-2}}
}{\sqrt{3\lambda _{\ell }-2}}
\end{array}
\right) =:\left(
\begin{array}{ccccc}
\mu_{1} & 0 & 0 & 0 & 0 \\
0 & \mu_{1} & 0 & 0 & 0 \\
0 & 0 & \mu_{3} & 0 & 0 \\
0 & 0 & 0 & \mu_{4} & 0 \\
0 & 0 & \mu_{2} & 0 & \mu_{5}
\end{array}
\right);
\end{equation*}
in the last expression, for notational simplicity we have omitted the
dependence of the $\mu_{i}$s on $\ell $. The matrix is block diagonal,
because under isotropy the gradient components are independent from the
Hessian when evaluated at the same point. We can hence define a $5$-dimensional
standard Gaussian vector $Y(x)=(Y_{1}(x),Y_{2}(x),Y_{3}(x),Y_{4}(x),Y_{5}(x))$
with independent components such that
\begin{eqnarray*}
	&&(\partial_{1}f_{\ell }(x),\partial_{2}f_{\ell }(x), \partial _{11}f_{\ell }(x),\partial _{12}f_{\ell }(x)-\cot \theta
	_{x}\partial _{2}f_{\ell }(x),\partial _{22}f_{\ell }(x)+\cot \theta
	_{x}\partial _{1}f_{\ell }(x)) \\
	&&=\Lambda _{\ell }Y(x)\\
	&&=\left( \mu_{1}Y_{1}(x),\mu_{1}Y_{2}(x),\mu
	_{3}Y_{3}(x),\mu_{4}Y_{4}(x),\mu_{5}Y_{5}(x)+\mu_{2}Y_{3}(x)\right).
\end{eqnarray*}
Hence
\begin{align*}
Y_{i}(x) &=\frac{\sqrt{2}}{\sqrt{\lambda _{\ell }}}\partial _{i;x}f_{\ell
}(x), \hspace{1cm} i=1,2, \\
Y_{3}(x) &=\frac{2\sqrt{2}}{\sqrt{\lambda _{\ell }}\sqrt{3\lambda _{\ell }-2
}}\partial _{11;x}f_{\ell }(x), \\
Y_{4}(x) &=\frac{2\sqrt{2}}{\sqrt{\lambda _{\ell }}\sqrt{\lambda _{\ell }-2}
}\partial _{21;x}, \\
Y_{5}(x) &=\frac{\sqrt{3\lambda _{\ell }-2}}{\lambda _{\ell }\sqrt{{\lambda
			_{\ell }-2}}}\partial _{22;x}f_{\ell }(x)-\frac{\lambda _{\ell }+2}{
	\lambda _{\ell }\sqrt{{\lambda _{\ell }-2}}\sqrt{3\lambda _{\ell }-2}}
\partial _{11;x}f_{\ell }(x).
\end{align*}
Note that asymptotically
\begin{equation*}
\mu_{1}\sim \frac{\ell }{\sqrt{2}}, \hspace{0.7cm} \mu_{2}\sim \frac{\ell ^{2}}{
	\sqrt{24}}, \hspace{0.7cm} \mu_{3}\sim \sqrt{\frac{3}{8}}\ell ^{2}, \hspace{0.7cm} \mu
_{4}\sim \frac{\ell ^{2}}{\sqrt{8}}, \hspace{0.7cm} \mu_{5}\sim \frac{\ell ^{2}}{\sqrt{3}},
\end{equation*}
where (as usual) $a_{\ell }\sim b_{\ell }$ means that the ratio between the
left- and right-hand side tends to unity as $\ell \rightarrow \infty$.
Thus we obtain
\begin{align*}
\mathcal{N}_{\ell }^{c}(I) &=\lim_{\varepsilon \to 0} \int_{\mathbb{S}^{2}}|\text{det}\nabla ^{2}f_{\ell
}(x)| \mathbb{I}_{\left\{ f_{\ell }(x)\in I\right\}} { \delta}_\varepsilon(\nabla f_{\ell }(x))dx \\
& =\lim_{\varepsilon \to 0} \lambda_{\ell }^{2}  \int_{\mathbb{S}^{2}}\left\vert
\frac{\mu_{3}\mu_{5}}{\lambda _{\ell }^{2}}Y_{3}(x)Y_{5}(x)+\frac{\mu_{2}\mu_{3}}{\lambda _{\ell }^{2}}Y_{3}^{2}(x)-\frac{\mu_{4}^{2}}{
	\lambda _{\ell }^{2}}Y_{4}^{2}(x)\right\vert  \mathbb{I}_{\{\frac{\mu_2+\mu_3}{\lambda_\ell }Y_{3}+\frac{\mu_5}{\lambda_\ell}Y_{5}\in I\}}  \\
	& \hspace{1.6cm} \times {\delta}_\varepsilon (\mu_1 Y_{1}(x),\mu_1 Y_{2}(x))dx,
\end{align*}
where
\begin{equation*}
\frac{\mu_{3}\mu_{5}}{\lambda _{\ell }^{2}}\sim \frac{1}{\sqrt{8}}, \hspace{0.7cm}
\frac{\mu_{2}\mu_{3}}{\lambda _{\ell }^{2}}\sim \frac{1}{8}, \hspace{0.7cm}
\frac{\mu_{4}^{2}}{\lambda _{\ell }^{2}}\sim \frac{1}{8}.
\end{equation*}
Since the $q$th order chaos is the space generated by the $L^{2}$-completion of
linear combinations of the form $H_{q_{1}}(Y_{1})\cdots H_{q_{5}}(Y_{5})$, with $
q_{1}+q_{2}+ \cdots +q_{5}=q$, it
is the linear span of cross-product of Hermite polynomials computed in the
independent random variables $Y_{i},$ $i=1, 2,\dots 5$,  which generate the
gradient and Hessian of $f_{\ell }$. In particular, the second order chaos
can be written in the following form
\begin{equation*}
{\mathcal{N}}_{\ell}^{c}(I)[2]= \lambda_{\ell}  \left[
\sum_{i<j}H_{ij}(I)\int_{\mathbb{S}^{2}}Y_{i}(x)Y_{j}(x)dx+\frac{1}{2}%
\sum_{i=1}^{5}K_{i}(I)\int_{\mathbb{S}^{2}}H_{2}(Y_{i}(x))dx\right],
\end{equation*}%
where
\begin{align*}
& H_{ij}(I)=\lim_{\varepsilon \rightarrow 0} \lambda_{\ell} \mathbb{E}\left[
\left\vert \frac{\mu_3\mu_5}{\lambda_\ell ^{2}}Y_{3}Y_{5}+\frac{%
	\mu_2\mu_3}{\lambda_\ell ^{2}}Y_{3}^{2}-\frac{\mu_4^{2}}{%
	\lambda_\ell ^{2}}Y_{4}^{2}\right\vert Y_{i}Y_{j}\; \id_{\{\frac{\mu_{2}+\mu_3}{\lambda_\ell }Y_{3}+\frac{\mu_5}{\lambda_\ell }Y_{5}\in
	I\}}\delta _{\varepsilon }(\mu_1 Y_{1}, \mu_1 Y_{2})\right] ,  \\ 
& K_{i}(I)=\lim_{\varepsilon \rightarrow 0} \lambda_{\ell} \mathbb{E}\left[ \left\vert
\frac{\mu_3\mu_5}{\lambda_\ell ^{2}}Y_{3}Y_{5}+\frac{\mu_{2}\mu_3}{\lambda_\ell ^{2}}Y_{3}^{2}-\frac{\mu_4^{2}}{\lambda_\ell ^{2}%
}Y_{4}^{2}\right\vert H_{2}(Y_{i})\; \id_{\{\frac{\mu_2+\mu_{3}}{\lambda_\ell }Y_{3}+\frac{\mu_5}{\lambda_\ell }Y_{5}\in I\}} \delta
_{\varepsilon }(\mu_1 Y_{1},\mu_1 Y_{2})\right];  
\end{align*}
the projection coefficients $H_{ij}(I)$ and $K_{i}(I)$ are constant with respect to $\ell$. The fourth order chaos is 
\begin{align} \label{4chaos}
\mathcal{N}_{\ell }^{c}(I)[4]
&=\lambda_{\ell } \left[ \frac{1}{2!2!}
\sum_{i=2}^{5}\sum_{j=1}^{i-1}h_{ij}(I) \int_{\mathbb{S}
	^{2}}H_{2}(Y_{i}(x))H_{2}(Y_{j}(x))dx+\frac{1}{4!}\sum_{i=1}^{5}k_{i}(I)
\int_{\mathbb{S}^{2}}H_{4}(Y_{i}(x))dx \right. \\
&\;\; +\left. \frac{1}{3!} \sum_{\stackrel{i,j=1}{i \ne j}}^5 g_{ij}(I) \int_{\mathbb{S}^{2}} H_3(Y_i(x)) H_1(Y_j(x)) d x \right. \nonumber \\ 
&+ \left. \frac{1}{2} \sum_{\stackrel{i,j,k=1}{i \ne j \ne k}}^5 p_{ijk}(I) \int_{\mathbb{S}^{2}} H_2(Y_i(x)) H_1(Y_j(x))H_1(Y_k(x)) d x   \right. \nonumber \\
&\;\; +\left. \sum_{\stackrel{i,j,k,l=1}{i \ne j \ne k \ne l}}^5 q_{ijkl}(I) \int_{\mathbb{S}^{2}} H_1(Y_i(x)) H_1(Y_j(x)) H_1(Y_k(x))  H_1(Y_l(x)) d x \right], \nonumber
\end{align}
where
\begin{align*}
h_{ij}(I)&=\lim_{\varepsilon\to 0}h_{ij}^\varepsilon(I)\\&=\lim_{\varepsilon \rightarrow 0} \lambda_{\ell } \; \mathbb{E}\left[ \left\vert
\frac{\mu_{3}\mu_{5}}{\lambda _{\ell }^{2}}Y_{3}Y_{5}+\frac{\mu
	_{2}\mu_{3}}{\lambda _{\ell }^{2}}Y_{3}^{2}-\frac{\mu_{4}^{2}}{\lambda
	_{\ell }^{2}}Y_{4}^{2}\right\vert  \mathbb{I}_{\{\frac{\mu_2+\mu
		_{3}}{\lambda_\ell }Y_{3}+\frac{\mu_5}{\lambda_\ell }Y_{5}\in I\}} {\delta}_{\varepsilon
}(\mu_1 Y_{1}, \mu_1Y_{2})H_{2}(Y_{i})H_{2}(Y_{j})\right], \\
k_{i}(I)&=\lim_{\varepsilon\to 0}k_{i}^\varepsilon(I)\\&
=\lim_{\varepsilon \rightarrow 0} \lambda_{\ell }\;  \mathbb{E}\left[ \left\vert
\frac{\mu_{3}\mu_{5}}{\lambda _{\ell }^{2}}Y_{3}Y_{5}+\frac{\mu
	_{2}\mu_{3}}{\lambda _{\ell }^{2}}Y_{3}^{2}-\frac{\mu_{4}^{2}}{\lambda
	_{\ell }^{2}}Y_{4}^{2}\right\vert  \mathbb{I}_{\{\frac{\mu_{2}+\mu_3}{\lambda_\ell }Y_{3}+\frac{\mu_5}{\lambda_\ell }Y_{5}\in
	I\}} { \delta}_{\varepsilon
}(\mu_1 Y_{1}, \mu_1 Y_{2})H_{4}(Y_{i})\right],\\
g_{ij}(I)&=\lim_{\varepsilon\to 0}g_{ij}^\varepsilon(I)\\&=\lim_{\varepsilon \rightarrow 0} \lambda_{\ell } \; \mathbb{E}\Big[ \left\vert
\frac{\mu_{3}\mu_{5}}{\lambda _{\ell }^{2}}Y_{3}Y_{5}+\frac{\mu
	_{2}\mu_{3}}{\lambda _{\ell }^{2}}Y_{3}^{2}-\frac{\mu_{4}^{2}}{\lambda
	_{\ell }^{2}}Y_{4}^{2}\right\vert \mathbb{I}_{\{\frac{\mu_{2}+\mu_3}{\lambda_\ell }Y_{3}+\frac{\mu_5}{\lambda_\ell }Y_{5}\in
	I\}} {\delta}_{\varepsilon
}(\mu_1 Y_{1}, \mu_1Y_{2}) \\
&  \hspace{1cm}\times H_3(Y_i(x)) H_1(Y_j(x)) \Big], \\
p_{ijk}(I)&=\lim_{\varepsilon\to} p_{ijk}^\varepsilon(I)\\&= \lim_{\varepsilon \rightarrow 0} \lambda_{\ell } \; \mathbb{E}\Big[ \left\vert
\frac{\mu_{3}\mu_{5}}{\lambda _{\ell }^{2}}Y_{3}Y_{5}+\frac{\mu
	_{2}\mu_{3}}{\lambda _{\ell }^{2}}Y_{3}^{2}-\frac{\mu_{4}^{2}}{\lambda
	_{\ell }^{2}}Y_{4}^{2}\right\vert \mathbb{I}_{\{\frac{\mu_{2}+\mu_3}{\lambda_\ell }Y_{3}+\frac{\mu_5}{\lambda_\ell }Y_{5}\in
	I\}} { \delta}_{\varepsilon} (\mu_1 Y_{1}, \mu_1Y_{2}) \\
&  \hspace{1cm}\times H_2(Y_i(x)) H_1(Y_j(x))H_1(Y_k(x))  \Big],\\
q_{ijkl}(I)&=\lim_{\varepsilon\to 0}q_{ijkl}^\varepsilon(I)\\&=\lim_{\varepsilon \rightarrow 0} \lambda_{\ell } \; \mathbb{E}\Big[ \left\vert
\frac{\mu_{3}\mu_{5}}{\lambda _{\ell }^{2}}Y_{3}Y_{5}+\frac{\mu
	_{2}\mu_{3}}{\lambda _{\ell }^{2}}Y_{3}^{2}-\frac{\mu_{4}^{2}}{\lambda
	_{\ell }^{2}}Y_{4}^{2}\right \vert \mathbb{I}_{\{\frac{\mu_{2}+\mu_3}{\lambda_\ell }Y_{3}+\frac{\mu_5}{\lambda_\ell }Y_{5}\in I \}} { \delta}_{\varepsilon }(\mu_1 Y_{1}, \mu_1Y_{2})  \\
&  \hspace{1cm}\times H_1(Y_i(x)) H_1(Y_j(x)) H_1(Y_k(x))  H_1(Y_l(x))\Big].
\end{align*}
\noindent The projection coefficients $k_{i}(I)$, $h_{ij}(I)$, $g_{ij}(I)$, $p_{ijk}(I)$, and $q_{ijkl}(I)$ do not depend on $\ell$. In \cite{CM2018b} it is proved that  
\begin{proposition}\label{values} For $I \subset \mathbb{R}$ and such that $\nu^c(I) \ne 0$, as $\ell \to \infty$, 
	\begin{align*}
	{\mathcal{N}}_{\ell}^{c}(I) -\mathbb{E}[{\mathcal{N}}_{\ell}^{c}(I)] &={\mathcal{N}}_{\ell}^{c}(I)[2] + R_{\ell}(I)
	\end{align*}
	where 
	\begin{align*}
	{\mathcal{N}}_{\ell}^{c}(I)[2] &= \frac{\lambda_{\ell}}{2} \nu^c(I) \frac{1}{2 \pi} \int_{S^2} H_2(f_{\ell}(x)) d x
	\end{align*}
	and 
	$${\rm Var}( {\mathcal{N}}_{\ell}^{c}(I)[2] )=\ell^3 [\nu^c(I)]^2+o(\ell^3), \hspace{1cm} \mathbb{E}[R^2_{\ell}(I)]=o(\ell^3),$$
	uniformly over $I$. 
\end{proposition}
Proposition \ref{values} says that the high frequency behaviour of the number of critical points is
dominated by a single term, proportional to the second-order Wiener chaos
projection and the second-order Wiener chaos projection
admits a simple expression in terms of the integral of $H_{2}(f_{\ell }(x))$
over $\mathbb{S}^{2}$. Recalling the definition of the random sequence
\begin{equation*}
\mathcal{F}_{\ell }=-\frac{\lambda _{\ell }}{2^{3}3^{2}\sqrt{3}\pi }\int_{
	\mathbb{S}^{2}}H_{4}(f_{\ell }(x))dx,
\end{equation*}
for which it is readily seen that
\begin{equation*}
\mathbb{E}\left[ \mathcal{F}_{\ell }\right] =0, \;\;\;\; \lim_{\ell
	\rightarrow \infty }\frac{\text{Var}(\mathcal{N}_{\ell }^{c})}{\text{Var}(\mathcal{F}
	_{\ell })}=1,
\end{equation*} 
in \cite{CM2019} it is proved that 
\begin{proposition} \label{points}
	As $\ell \rightarrow \infty $
	\begin{equation*}
	{\rm Corr}(\mathcal{N}_{\ell }^{c},\mathcal{F}_{\ell }) \rightarrow 1,
	\end{equation*}
	and hence
	\begin{equation*}
	\mathcal{N}_{\ell }^{c}-\mathbb{E}\left[ \mathcal{N}_{\ell }^{c}\right]
	=\mathcal{F}_{\ell}+o_{\mathbb{P}}(\sqrt{\ell^2 \log \ell}).
	\end{equation*}
\end{proposition}
As a consequence, the total number of critical points is
fully correlated in the limit, for $\ell \rightarrow \infty$, with $h_{\ell ;4}$. The
limiting distribution of $h_{\ell ;4}$ is studied
in \cite{MW2014}, where it is shown that a quantitative version of the
Central Limit Theorem holds for $h_{\ell ;4}$. Proposition \ref{values} also implies that the asymptotic behaviour of the total number of critical points is dominated by its projection on the fourth-order chaos term and that the projection on the fourth-order chaos can be expressed simply in terms
	of the fourth-order Hermite polynomial, evaluated on the eigenfunctions $\left\{
	f_{\ell }\right\}$, without the need to compute Hermite polynomials
	evaluated on the first and second derivatives of $\left\{f_{\ell }\right\}$, despite the fact that the
	latter do appear in the Kac-Rice formula and they are not negligible in
	terms of asymptotic variance.
A consequence is that, as $\ell \rightarrow \infty$,
\begin{equation*}
\text{Var}\left( \mathcal{N}_{\ell }^{c}[4]\right)=\frac{\ell ^{2}\log \ell }{
	3^{3}\pi ^{2}}+O(\ell ^{2}),
\end{equation*}
so that
\begin{equation*}
\lim_{\ell \rightarrow \infty }\frac{\text{Var}\left( \mathcal{N}_{\ell
	}^{c}[4]\right) }{\text{Var}\left( \mathcal{N}_{\ell }^{c}\right) }=1.
\end{equation*}
Note that by orthogonality we have
\begin{equation*}
\text{Var}\left( \mathcal{N}_{\ell }^{c}\right) =\sum_{q=0}^{\infty }\text{Var}\left(
\mathcal{N}_{\ell }^{c}[q]\right) =\text{Var}\left(
\mathcal{N}_{\ell }^{c}[4]\right) +\sum_{k=1}^{\infty }\mathcal{N}_{\ell }^{c}[4+2k],
\end{equation*}
where the odd terms in the expansion vanish by symmetry arguments,
$\rm{Var}\left( \mathcal{N}_{\ell }^{c}[0]\right) =0$ is obvious, and
$\rm{Var}\left( \mathcal{N}_{\ell }^{c}[2]\right) =0$ is shown in
\cite{CM2018b}. Hence we have the asymptotic relation
\begin{equation*}
\sum_{k=1}^{\infty }\mathcal{N}_{\ell }^{c}[4+2k]=o(\ell^{2}\log \ell ).
\end{equation*}

\section{Proof of Theorem \ref{main theorem}}

In this section we give the proof of our main result. Let us consider $I_1, I_2 \subset \mathbb{R}$; assume first $I_i \ne \mathbb{R}$ and such that $\nu^c(I_i)\ne 0$, $i=1,2$. Thanks to Proposition \ref{values} we have 
\begin{align} \label{15:43}
{\mathcal{N}}_{\ell}^{c}(I_i) -\mathbb{E}[{\mathcal{N}}_{\ell}^{c}(I_i)] &={\mathcal{N}}_{\ell}^{c}(I_i)[2] + R_{\ell}(I_i)=\mathcal{S}_{\ell}(I_i) + R_{\ell}(I_i),
\end{align}
where $\mathbb{E}[R^2_{\ell}(I)]=o(\ell^3)$, as $\ell \to \infty$.  Moreover, recalling that 
$$\text{Var}( {\mathcal{N}}_{\ell}^{c}(I_i)[2] )=\ell^3 [\nu^c(I_i)]^2+o(\ell^3)$$ 
and, see \cite{CMW14}, 
$${\rm Var}(\mathcal{N}^c_\ell(I_i))=\ell^3 [\nu^c(I_i)]^2+O(\ell^{5/2}),$$ 
we find that, as $\ell \to \infty$, 
$${\rm Corr}(\mathcal{N}_\ell^c(I_1),\mathcal{N}_\ell^c(I_2)) \to 1.$$
We consider now the total number of critical points, in view of Proposition \ref{points} we have 
\begin{equation} \label{15:44}
\mathcal{N}_{\ell }^{c}-\mathbb{E}\left[ \mathcal{N}_{\ell }^{c}\right]
=\mathcal{N}_{\ell }^{c}[4]+ R_{\ell} =\mathcal{F}_{\ell}+R_{\ell}, 
\end{equation}
where 
$R_{\ell}=o_{\mathbb{P}}(\sqrt{\ell^{2}\log \ell })$ as $\ell \to \infty$. From \eqref{15:43} and \eqref{15:44}  we have that for $I \ne \mathbb{R}$ and $I$ such that $\nu^c(I_i)\ne 0$ 
$${\rm Corr}(\mathcal{N}_\ell^c(I),\mathcal{N}_\ell^c) \to 0.$$
Now we observe that, taking 
$$Z=||f_{\ell}||^2_{L^2(\mathbb{S}^2)},$$ 
we have 
$$\mathcal{N}_\ell^c(I)^*=\sum_{q=3}^{\infty} \mathcal{N}_\ell^c(I)[q].$$
We state now the following results:  
\begin{proposition}\label{VarNu2} As $\ell \to \infty$, 
	\begin{equation*}\label{Var critical values}
	{\rm Var}(\mathcal{N}_\ell^c(I))=\frac{[ 5\mathcal{I}_0(I)- \mathcal{I}_2(I) ]^2}{2^{4}} \ell^3+ \frac{[51 \mathcal{I}_0(I)-2 \cdot 11\mathcal{I}_2(I)+\mathcal{I}_4(I)]^2}{2^{6}\pi^2}\ell^2 \log \ell+O(\ell^2),
	\end{equation*}
	where
	\begin{align*}
	\mathcal{I}_0(I)&=\sqrt{\frac{2}{\pi}} \int_I  [2 e^{-t^2} + t^2-1] e^{- \frac{t^2}{2}} dt, \\
	\mathcal{I}_2(I)&= \sqrt{\frac{2}{\pi}} \int_I  [-4+t^2+t^4+ e^{-t^2} 2(4+ 3 t^2)] e^{- \frac{t^2}{2}} dt,  \\
	\mathcal{I}_4(I)&= \sqrt{\frac{2}{\pi}} \int_I  [(72+96 t^2+38 t^4) e^{-t^2} -36 -12 t^2 +11 t^4+t^6] e^{-\frac{t^2}{2}} dt.   
	\end{align*}
\end{proposition}
The proof of Proposition \ref{VarNu2} is postponed to Section \ref{19:20}. The next corollary follows immediately from Proposition \ref{VarNu2} and 
$$\text{Var}( {\mathcal{N}}_{\ell}^{c}(I_i)[2] )=\ell^3 [\nu^c(I_i)]^2+o(\ell^3),$$
by observing that 
$$\nu^c(I)=\frac{ 5\mathcal{I}_0(I)- \mathcal{I}_2(I) }{2^{2}}.$$
\begin{corollary} \label{lem1} As $\ell \to \infty$
	$${\rm Var} \Big(\sum_{q=3}^{\infty} \mathcal{N}_\ell^c(I)[q] \Big)=\frac{[51 \mathcal{I}_0(I)-2 \cdot 11\mathcal{I}_2(I)+\mathcal{I}_4(I)]^2}{2^{6}\pi^2}\ell^2 \log \ell+O(\ell^2). $$
\end{corollary} 
Let us denote by 
$$\mathcal{F}_\ell(I)= \frac{51 \mathcal{I}_0(I)-2 \cdot 11\mathcal{I}_2(I)+\mathcal{I}_4(I)}{2^{3} \pi } \lambda_\ell \int_{\mathbb{S}^2}H_4(f_\ell(x)) \,dx.$$
Note that $\mathcal{F}_\ell=\mathcal{F}_\ell(\mathbb{R})$. We have the following proposition whose proof is given in Section \ref{prooflem2}.
\begin{proposition} \label{lem2}
	For all $I \subset \mathbb{R}$, as $\ell \to \infty$, 
	$${\rm Corr} (\sum_{q=3}^{\infty} \mathcal{N}_\ell^c(I)[q], \mathcal{F}_{\ell}(I)) \to 1.$$
\end{proposition}
The statement of Theorem \ref{main theorem} follows from Corollary \ref{lem1} and Proposition \ref{lem2}: for every choice of $I_1$ and $I_2$, we have that 
$${\rm Corr}_{||f_{\ell} ||^2_{L^2(\mathbb{S}^2)}}(\mathcal{N}_\ell^c(I_1),\mathcal{N}_\ell^c(I_2)) \to 1.$$

\section{Proof of Proposition \ref{VarNu2}} \label{19:20}

\subsection{Approximate Kac-Rice formula for counting critical points with value in $I \subseteq \R$} \label{K-RinI}
For counting the critical points with corresponding value lying in an interval $I$ in the real line, we define the two-point correlation function: for $x \ne \pm y$ 
\begin{equation*} 
\begin{split}
&K_{2,\ell}(x,y;t_{1},t_{2})\\&=
\mathbb{E}\left[  \left\vert \nabla ^{2}f_{\ell }(x)\right\vert
\cdot \left\vert \nabla ^{2}f_{\ell }(y)\right\vert \Big|  \nabla
f_{\ell }(x)=\nabla f_{\ell }(y)= 0,f_{\ell }(x)=t_{1},f_{\ell
}(y)=t_{2}\right] \cdot \varphi _{x,y, \ell}(t_{1},t_{2}, 0,0),
\end{split}
\end{equation*}
where $\varphi _{x,y, \ell}(t_{1},t_{2},0,0)$ denotes the density of the $6$-dimensional vector
$$
\left( f_{\ell }(x),f_{\ell }(y),\nabla f_{\ell }(x),\nabla
f_{\ell }(y)\right)
$$
in $f_{\ell }(x)=t_{1},f_{\ell }(y)=t_{2},\nabla f_{\ell }(x)=\nabla f_{\ell}(y)=0$. In  \cite{CMW14} the following {\it approximate} Kac-Rice formula is derived: there exists a constant $C>0$ sufficiently big such that 
\begin{align} \label{afkaoI}
{\rm Var} \left( {\cal N}^c_\ell(I) \right) =  \int_{\mathcal{W}} \iint_{I \times I} K_{2,\ell}(x,y; t_1, t_2)  \; d t_1 d t_2 d x d y - (\E [ {\cal N}^c_\ell(I)])^2 +  O(\ell^2),
\end{align}
where $\mathcal{W}$ is the union of all tuples of points belonging to Voronoi cells further apart than $C/\ell$, see \cite{CMW14} for a more formal definition. Now, exploiting isotropy, and observing that the level field $f_{\ell}$ is a linear combination of gradient and second order derivatives, 
we have, see \cite[Section 4.1.2]{CMW14}: 
\begin{equation*}
K_{2,\ell }(\phi ; t_{1},t_{2}) = \frac{\lambda_{\ell}^{4}}{8} \frac{1}{\pi^{2} \sqrt{(\lambda _{\ell }^{2}-4\alpha _{2,\ell }^{2}(\phi
		))(\lambda _{\ell }^{2}-4\alpha _{1,\ell }^{2}(\phi ))}} q (\mathbf{a}_{\ell}(\phi);t_1,t_2),
\end{equation*}
where $\phi=d(x,y)$ is the geodesic distance between the two points $x$ and $y$, and  
\begin{align*}
&q (\mathbf{a}_{\ell}(\phi);t_1,t_2)\\
&=\frac{1}{(2\pi )^{3}   \sqrt{\det (\Delta (\mathbf{a}_{\ell}(\phi)))}  }\iint_{\mathbb{R}^{2}\times \mathbb{R}%
	^{2}} \left|z_{1} \sqrt 8 t_1 -z_{1}^{2}-z_{2}^{2}\right|\cdot \left| w_{1} \sqrt 8 t_2 -w_{1}^{2}-w_{2}^{2}\right| \\
& \times  \exp \left\{-\frac{1}{2} 
(z_{1},z_{2},\sqrt 8
t_{1}-z_{1},w_{1},w_{2}, \sqrt 8 t_{2}-w_{1}) \Delta
(\mathbf{a}_{\ell}(\phi))^{-1}   (z_{1},z_{2},\sqrt 8
t_{1}-z_{1},w_{1},w_{2}, \sqrt 8 t_{2}-w_{1})^t   \right\}\\
& \times d z_1 dz_2 d w_1 dw_2.
\end{align*}
The matrix $\Delta (\mathbf{a}_{\ell}(\phi))$ is defined as 
$$\Delta (\mathbf{a}_{\ell}(\phi))=\frac{8}{\lambda^2_{\ell}} \Omega_{\ell}(\phi)$$
where $\Omega_{\ell}(\phi)$ is the conditional covariance matrix of the random vector 
$$(\nabla^2 f_{\ell}(x), \nabla^2 f_{\ell}(y) | \nabla f_{\ell}(x)=\nabla f_{\ell}(y)=0),$$   
 see \cite[Section 4.1.1]{CMW14}. 

\subsection{Taylor expansion} In \cite[Section 3.3]{CW15}, it is proved that as $\ell \to \infty$, $\text{Var}(\mathcal{N}^c_\ell(I))$ has the following leading terms
\begin{align*} 
\text{Var}(\mathcal{N}^c_\ell(I))&= \frac 1 8 \left[ 2 \ell^3 + \frac{2 \cdot 3^2}{ \pi^2} \ell^2 \log \ell \right]  \iint_{I \times I} q(\mathbf{0};t_1,t_2) dt_1 dt_2 \\
&+ \frac 1 8 \left[ - 16 \ell^3 - \frac{2^5 \cdot 3 }{\pi^2} \ell^2 \log \ell \right] \; \iint_{I \times I} \Big[\frac{\partial }{\partial a_{3} }  q (\mathbf{a};t_1,t_2)\Big]_{\mathbf{a} =\mathbf{0}} dt_1 dt_2   \nonumber \\
&+ \frac 1 8 \left[ 32 \ell^3 - \frac{ 2^6}{\pi^2} \ell^2 \log \ell \right] 
 \; \iint_{I \times I} \Big[ \frac{\partial^2 }{\partial a_{7} \partial a_{7} }  q (\mathbf{a};t_1,t_2)\Big]_{\mathbf{a} =\mathbf{0}} dt_1 dt_2 \\
& + \frac 1 8 \left[ \frac{ 3 \cdot 2^7}{\pi^2} \ell^2 \log \ell \right]  \; \iint_{I \times I} \Big[ \frac{\partial^2 }{\partial a_{3} \partial a_{3} }  q (\mathbf{a};t_1,t_2)\Big]_{\mathbf{a} =\mathbf{0}} dt_1 dt_2 \nonumber \\
&+ \frac 1 8 \left[ - \frac{ 2^9}{ \pi^2} \ell^2 \log \ell \right] \; \iint_{I \times I} \Big[\frac{\partial^3 }{\partial a_{3} \partial a_{7}  \partial a_{7} }  q (\mathbf{a};t_1,t_2)\Big]_{\mathbf{a} =\mathbf{0}} dt_1 dt_2 \\
&+ \frac 1 8 \left[ \frac{ 2^9}{ \pi^2} \ell^2 \log \ell \right] \; \iint_{I \times I} \Big[\frac{\partial^4 }{\partial a_{7}\partial a_{7} \partial a_{7}  \partial a_{7} }  q (\mathbf{a};t_1,t_2)\Big]_{\mathbf{a} =\mathbf{0}} dt_1 dt_2 +O(\ell^2).
\end{align*}
 
Let $Z=(Z_{1},Z_{2},Z_{3})$ be a centred jointly Gaussian random vector with covariance matrix
\begin{equation*}
\tilde{c}_{\infty }=\left(
\begin{array}{ccc}
3 & 0 & 1 \\
0 & 1 & 0 \\
1 & 0 & 3%
\end{array}%
\right), 
\end{equation*}%
we denote by $\phi _{Z_{1}+Z_{3}}$ the probability density function of $Z_{1}+Z_{3}$ and we define 
\begin{align*}
p_0(t)=\sqrt{8} \cdot \mathbb{E} [|Z_1 Z_3 - Z_2^2| \Big|  Z_1+Z_3=\sqrt{8} t] \cdot \phi_{Z_1+Z_3} (\sqrt{8} t).
\end{align*} 
An explicit computation, see \cite[Remark 4.1]{CMW14}, shows that 
\begin{align*}
\iint_{I \times I} q(\mathbf{0};t_1,t_2) dt_1 dt_2=\frac{1}{2^3} \left[ \int_I p_0(t) d t \right]^2,
\end{align*} 
where
\begin{align*}
p_0(t)= \sqrt{\frac{2}{\pi}} [2 e^{-t^2} + t^2-1] e^{- \frac{t^2}{2}}.
\end{align*} 
Now let ${\cal I}_{r}(I)=\int_{I} p_r(t) dt$, $r=0,2,4$, with  $p_r$, $r=2,4$, defined by 
\begin{align*}
p_2(t)&=\sqrt{8} \cdot \mathbb{E} [(3 t -\sqrt 2 Z_1)^2 |Z_1 Z_3 - Z_2^2| \Big|  Z_1+Z_3=\sqrt{8} t] \cdot \phi_{Z_1+Z_3} (\sqrt{8} t) \\
&= \sqrt{\frac{2}{\pi}} [-4+t^2+t^4+ e^{-t^2} 2(4+ 3 t^2)] e^{- \frac{t^2}{2}},\\
p_4(t)&=\sqrt{8} \cdot \mathbb{E} [(3 t -\sqrt 2 Z_1)^4 |Z_1 Z_3 - Z_2^2| \Big|  Z_1+Z_3=\sqrt{8} t] \cdot \phi_{Z_1+Z_3} (\sqrt{8} t) \\
&=  \sqrt{\frac{2}{\pi}} [(72+96 t^2+38 t^4) e^{-t^2} -36 -12 t^2 +11 t^4+t^6] e^{-\frac{t^2}{2}},
\end{align*}
using Leibnitz integral rule and some 
technical computations, we obtain the following: 
\begin{align*}
& \iint_{I \times I} \Big[\frac{\partial }{\partial a_{3} }  q (\mathbf{a};t_1,t_2)\Big]_{\mathbf{a} =\mathbf{0}} dt_1 dt_2 = -\frac{3}{2^6} {\cal I}^2_{0}(I)+\frac{1}{2^6} {\cal I}_{0}(I) {\cal I}_{2}(I), \\
&\iint_{I \times I} \Big[\frac{\partial^2 }{\partial a_{7} \partial a_{7} }  q (\mathbf{a};t_1,t_2)\Big]_{\mathbf{a} =\mathbf{0}} dt_1 dt_2  
=  \frac{1}{2^9} [3 {\cal I}_{0}(I) - {\cal I}_{2}(I)   ]^2,\\
&\iint_{I \times I} \Big[\frac{\partial^2 }{\partial a_{3} \partial a_{3} }  q (\mathbf{a};t_1,t_2)\Big]_{\mathbf{a} =\mathbf{0}} dt_1 dt_2
= \frac{1}{2^{11}} [    2^3  \cdot 3^2  {\cal I}^2_{0}(I)  - 2^4 \cdot 3  {\cal I}_{0}(I)  {\cal I}_{2}(I)  + 2 {\cal I}_{0}(I)   {\cal I}_{4}(I)+2  {\cal I}^2_{2}(I) ], \\
&\iint_{I \times I} \Big[\frac{\partial^3 }{\partial a_{3} \partial a_{7}  \partial a_{7} }  q (\mathbf{a};t_1,t_2)\Big]_{\mathbf{a} =\mathbf{0}} dt_1 dt_2 \\&= \frac{1}{2^{13}} [ -2 \cdot 3^4 {\cal I}_{0}(I)^2 +2 \cdot  3^4 {\cal I}_{0}(I) {\cal I}_{2}(I) - 2 \cdot 3 {\cal I}_{0}(I) {\cal I}_{4}(I)  + 2 {\cal I}_{4}(I) {\cal I}_{2}(I)  -2^2 \cdot 3^2 {\cal I}^2_{2}(I)  ], \\
&\iint_{I \times I} \Big[\frac{\partial^4 }{\partial a_{7}\partial a_{7} \partial a_{7}  \partial a_{7} }  q (\mathbf{a};t_1,t_2)\Big]_{\mathbf{a} =\mathbf{0}} dt_1 dt_2
=\frac{1}{2^{15}} [3^3 {\cal I}_{0}(I) - 2 \cdot 3^2 {\cal I}_{2}(I)+{\cal I}_{4}(I)]^2.
\end{align*}
Therefore we find the following analytic expression for the variance:
\begin{align*}
&\text{Var}({\cal N}^c_\ell(I))=\left[ 2 \ell^3 + \frac{2 \cdot 3^2}{\pi^2} \ell^2 \log \ell \right]   \frac{1}{2^3} {\cal I}^2_{0}(I) \\
&+\left[ - 16  \ell^3 - \frac{2^5 \cdot 3 }{\pi^2} \ell^2 \log \ell \right]   \;   \frac{1}{2^6} [-3 {\cal I}^2_{0}(I)+\ {\cal I}_{0}(I) {\cal I}_{2}(I) ]+\left[ 32 \ell^3 - \frac{ 2^6}{\pi^2} \ell^2 \log \ell  \right]   \frac{1}{2^9} [3 {\cal I}_{0}(I) - {\cal I}_{2}(I)   ]^2 \\
&+\left[  \frac{ 3 \cdot 2^7}{\pi^2} \ell^2 \log \ell  \right]   \frac{1}{2^{11}} [    2^3  \cdot 3^2  {\cal I}^2_{0}(I)  - 2^4 \cdot 3  {\cal I}_{0}(I)  {\cal I}_{2}(I)  + 2 {\cal I}_{0}(I)   {\cal I}_{4}(I)+2  {\cal I}^2_{2}(I)      ] \\
&+\left[ - \frac{  2^9}{ \pi^2} \ell^2 \log \ell  \right]  \frac{1}{2^{13}} [ -2 \cdot 3^4 {\cal I}_{0}(I)^2 +2 \cdot  3^4 {\cal I}_{0}(I) {\cal I}_{2}(I) - 2 \cdot 3 {\cal I}_{0}(I) {\cal I}_{4}(I)  + 2 {\cal I}_{4}(I) {\cal I}_{2}(I)  -2^2 \cdot 3^2 {\cal I}^2_{2}(I)  ] \\
&+\left[  \frac{  2^9}{ \pi^2} \ell^2 \log \ell  \right] \frac{1}{2^{15}} [3^3 {\cal I}_{0}(I) - 2 \cdot 3^2 {\cal I}_{2}(I)+{\cal I}_{4}(I)]^2+O(\ell^2),
\end{align*}
that can be rewritten as: 
\begin{align*}
\text{Var}({\cal N}^c_\ell(I))= \frac{1}{2^4}[ 5 {\cal I}_{0}(I) - {\cal I}_{2}(I)]^2\;  \ell^3+  \frac{1}{\pi^2 2^{6}}  [51 {\cal I}_{0}(I) -  2\cdot 11\; {\cal I}_{2}(I)+{\cal I}_{4}(I)  ]^2\; \ell^2 \log \ell +O(\ell^2),
\end{align*}
where
\begin{align*}
&5 {\cal I}_{0}(I) - {\cal I}_{2}(I) = \sqrt{\frac{2}{\pi}} \int_I e^{-\frac{t^2}{2} }  [- 1+4 t^2 -t^4 +e^{-t^2}(2-6 t^2) ]  d t , \\
&51 {\cal I}_{0}(I) -  2\cdot 11\; {\cal I}_{2}(I)+{\cal I}_{4}(I)   = \sqrt{\frac{2}{\pi}} \int_I  e^{-\frac{t^2}{2}}  [(-2 -36 t^2+38 t^4) e^{-t^2} +1+17 t^2-11 t^4+t^6] d t,
\end{align*}
that is the statement of Proposition \ref{VarNu2}.

\section{Proof of Proposition \ref{lem2}}\label{prooflem2}

In order to prove Proposition \ref{lem2} we first give the following results, whose proofs are given in Section \ref{EvaCoeff} and Section \ref{oddt} respectively.
\begin{proposition}
	\label{coefficients} We have that
	\begin{align*}
	k_{2}(I)&=\frac{3}{8\pi } \cal{I}_0(I), \\
	k_{5}(I)&=\frac{3}{8\pi }\cal{I}_0(I)-\frac{2^3}{2^5\pi} \cal{I}_2(I)+ \frac{2^6}{3^2 2^9 \pi} \cal{I}_4(I),\\
	h_{25}(I)&=\frac{1}{2^3 \pi} \cal{I}_0(I)-\frac{2^3}{2^6 3 \pi} \cal{I}_2(I).
	\end{align*}
\end{proposition}

\begin{proposition} \label{17:37}
	The projection coefficients $g_{ij}(I)$, $p_{ijk}(I)$ and $q_{ijkl}(I)$ are such that
	\begin{itemize}
		\item for $i, j \ne 3, 5$, we have  $g_{ij}(I)=0$,
		\item for $j, k \ne 1, 2, 4$, we have $p_{ijk}(I)=0$,
		\item we have $q_{ijkl}(I)=0$.
	\end{itemize}
\end{proposition}

We also recall in the two lemmas below the results given in \cite[Lemmas 5.1, 5.2, 5.3, 5.4, 5.6]{CM2019}.

\begin{lemma}\label{Lem5.1-5.2-5.3}
	As $\ell \to \infty$, 
	\begin{equation*}
	\int_{0}^{\pi /2}\mathbb{E}\left[ H_{4}(Y_{2}(\bar{x}))H_{4}(f_{\ell
	}(y(\phi )))\right] \sin \phi d\phi =4!\frac{2 \cdot 3}{\pi ^{2}}\frac{\log \ell }{\ell ^{2}}+O(\ell ^{-2}),
	\end{equation*}
	\begin{equation*}
	\int_{0}^{\pi /2}\mathbb{E}\left[ H_{4}(Y_{5}(\bar{x}))H_{4}(f_{\ell
	}(y(\phi )))\right] \sin \phi d\phi =4!\frac{3^3}{2\pi ^{2}}\frac{\log
		\ell }{\ell ^{2}}+O(\ell ^{-2}),
	\end{equation*}
	\begin{equation*}
	\int_{0}^{\pi /2}\mathbb{E}\left[ H_{2}(Y_{2}(\bar{x}))H_{2}(Y_{5}(\bar{x}
	))H_{4}(f_{\ell }(y(\phi )))\right] \sin \phi d\phi =4!\frac{3}{\pi ^{2}}
	\frac{\log \ell }{\ell ^{2}}+O(\ell ^{-2}).
	\end{equation*}
	
	
\end{lemma}

\begin{lemma}\label{lem5.4,5.5,5.6}
As $\ell \to \infty$, 

\begin{align*}
\int_{0}^{\pi / 2} \mathbb{E}[H_{3}(Y_{a}(\bar{x}))H_{1}(Y_{b}(\bar{x})) H_4(f_{\ell}(y(\phi)))] \sin \phi d \phi= O(\ell^{-2}),
\end{align*}
for $i=1, 2$,
\begin{align*}
\int_{0}^{\pi / 2} \mathbb{E}[H_{2}(Y_{i}(\bar{x})) H_{1}(Y_{3}(\bar{x})) H_{1}(Y_{5}(\bar{x}))  H_4(f_{\ell}(y(\phi)))] \sin \phi d \phi=0,
\end{align*}
and 
\begin{align*}
\int_{0}^{\pi / 2} \mathbb{E}[H_{2}(Y_{4}(\bar{x})) H_{1}(Y_{3}(\bar{x}))  H_{1}(Y_{5}(\bar{x})) H_4(f_{\ell}(y(\phi)))] \sin \phi d \phi=0.
\end{align*}
\end{lemma}

Relying on the results just mentioned we can now prove Proposition \ref{lem2}:  by continuity of the inner product in $L^{2}(\Omega )$, we write
\begin{align*}
\text{Cov}(\mathcal{N}_{\ell }^{c}(I),h_{\ell ;4}) &=\lim_{\varepsilon \rightarrow
	0}\text{Cov}(\mathcal{N}_{\ell ,\varepsilon }^{c}(I),h_{\ell ;4})  \\
&=\lim_{\varepsilon \rightarrow 0}\mathbb{E}\left[ \int_{\mathbb{S}^{2}}|
\text{det} \nabla ^{2}f_{\ell }(x)| \mathbb{I}_{\left\{ f_{\ell }(x)\in I\right\}} { \delta}_{\varepsilon }(\nabla f_{\ell
}(x))dx\int_{\mathbb{S}^{2}}H_{4}(f_{\ell }(y))dy\right] \\
&=\lim_{\varepsilon \rightarrow 0}\mathbb{E}\left[ \int_{\mathbb{S}
	^{2}}\sum_{q=0}^{\infty }\psi _{\ell }^{\varepsilon }(I,x ;q)dx\int_{\mathbb{S}
	^{2}}H_{4}(f_{\ell }(y))dy\right].
\end{align*}
Now note that both $\psi_{\ell }^{\varepsilon }(I,x;q)$ and $H_{4}(f_{\ell
}(y)) $ are isotropic processes on $\mathbb{S}^{2}$, hence we have
\begin{align*}
\mathbb{E}\left[ \int_{\mathbb{S}^{2}}\sum_{q=0}^{\infty }\psi _{\ell
}^{\varepsilon }(I,x;q)dx\int_{\mathbb{S}^{2}}H_{4}(f_{\ell }(y))dy\right]
&=\mathbb{E}\left[ \int_{\mathbb{S}^{2}}\lim_{Q\rightarrow \infty
}\sum_{q=0}^{Q}\psi _{\ell }^{\varepsilon }(I,x;q)dx\int_{\mathbb{S}
	^{2}}H_{4}(f_{\ell }(y))dy\right] \\
&=\lim_{Q\rightarrow \infty }\mathbb{E}\left[ \int_{\mathbb{S}
	^{2}}\sum_{q=0}^{Q}\psi _{\ell }^{\varepsilon }(I,x;q)dx\int_{\mathbb{S}
	^{2}}H_{4}(f_{\ell }(y))dy\right]
\end{align*}
by continuity of covariances. Moreover because all integrands are finite-order polynomials we have
\begin{align*}
\lim_{Q\rightarrow \infty }\mathbb{E}&\left[ \int_{\mathbb{S}
	^{2}}\sum_{q=0}^{Q}\psi _{\ell }^{\varepsilon }(I,x;q)dx\int_{\mathbb{S}
	^{2}}H_{4}(f_{\ell }(y))dy\right]
=\lim_{Q\rightarrow \infty }\sum_{q=0}^{Q}\int_{\mathbb{S}^{2}}\int_{
	\mathbb{S}^{2}}\mathbb{E}\left[ \psi_{\ell }^{\varepsilon}(I,x;q)H_{4}(f_{\ell }(y))\right] dxdy  \\
&=\int_{\mathbb{S}^{2}}\int_{\mathbb{S}^{2}}\mathbb{E}\left[ \psi_{\ell
}^{\varepsilon }(I,x;4)H_{4}(f_{\ell }(y))\right] dxdy =16\pi ^{2}\int_{0}^{\pi /2}\mathbb{E}\left[ \psi_{\ell }^{\varepsilon }(I,\overline{x};4)H_{4}(f_{\ell }(y(\phi )))\right] \sin \phi d\phi,
\end{align*}
where in the last steps we used orthogonality of Wiener chaoses and isotropy;
we take $\overline{x}=(\frac{\pi }{2},0)$ and $y(\phi )=(\frac{\pi }{2},\phi )$.
This allows us to perform our argument on the {\it equator}, where $\theta $ is fixed to $\pi/2$.
Note that
\begin{align*}
\psi_{\ell }^{\varepsilon }(I,\overline{x};4)&=\lambda_{\ell} \Big[\frac{1}{2!2!}\sum_{i=2}^{5}\sum_{j=1}^{i-1}h_{ij}^{
	\varepsilon }(I)H_{2}(Y_{i}(\bar{x}))H_{2}(Y_{j}(\bar{x}))+\frac{1}{4!}
\sum_{i=1}^{5}k_{i}^{\varepsilon }(I)H_{4}(Y_{i}(\bar{x}))\\
&\;\;+ \frac{1}{3!} \sum_{\stackrel{i,j=1}{i \ne j}}^5 g_{ij}^{\varepsilon } (I) H_3(Y_i(\bar{x})) H_1(Y_j(\bar{x}))  + \frac{1}{2} \sum_{\stackrel{i,j,k=1}{i \ne j \ne k}}^5 q_{ijk}^{\varepsilon } (I) H_2(Y_i(\bar{x})) H_1(Y_j(\bar{x}))H_1(Y_k(\bar{x}))   \\
&\;\; + \sum_{\stackrel{i,j,k,l=1}{i \ne j \ne k \ne l}}^5 p^{\varepsilon }_{ijkl} (I) H_1(Y_i(\bar{x})) H_1(Y_j(\bar{x})) H_1(Y_k(\bar{x}))  H_1(Y_l(\bar{x})) \Big],
\end{align*}
and hence
\begin{align*}
&\text{Cov}(\mathcal{N}_{\ell }^{c}(I),h_{\ell ;4}) \\
&=16\pi ^{2}\lim_{\varepsilon
	\rightarrow 0}\int_{0}^{\pi /2}\mathbb{E}\left[ \psi_{\ell }^{\varepsilon }(I, \overline{x};4)H_{4}(f_{\ell }(y(\phi )))\right] \sin \phi d\phi \\
&=16\pi^{2} \lambda_{\ell}  \frac{1}{2! 2!}
\sum_{i=2}^{5}\sum_{j=1}^{i-1} \{\lim_{\varepsilon \rightarrow
	0}h_{ij}^{\varepsilon } (I)\} \int_{0}^{\pi /2}\mathbb{E}\left[ H_{2}(Y_{i}(
\bar{x}))H_{2}(Y_{j}(\bar{x}))H_{4}(f_{\ell }(y(\phi )))\right] \sin \phi
d\phi \\
&\;\;+16\pi^{2} \lambda_{\ell} \frac{1}{4!}
\sum_{i=1}^{5} \{\lim_{\varepsilon \rightarrow 0}k_{i}^{\varepsilon }(I)
\} \int_{0}^{\pi /2}\mathbb{E}\left[ H_{4}(Y_{i}(\bar{x}))H_{4}(f_{\ell
}(y(\phi )))\right] \sin \phi d\phi\\
&\;\;+16\pi^{2} \lambda_{\ell} \frac{1}{3!}
\sum_{i \ne j} \{\lim_{\varepsilon \rightarrow 0}g_{ij}^{\varepsilon }(I)
\} \int_{0}^{\pi /2} \mathbb{E}\left[ H_{3}(Y_{i}(\bar{x})) H_{1}(Y_{j}(\bar{x})) H_{4}(f_{\ell
}(y(\phi )))\right] \sin \phi d\phi\\
&\;\;+16\pi^{2} \lambda_{\ell} \frac{1}{2}
\sum_{i \ne j \ne k} \{\lim_{\varepsilon \rightarrow 0} p_{ijk}^{\varepsilon }(I)
\} \int_{0}^{\pi /2}\mathbb{E}\left[ H_{2}(Y_{i}(\bar{x})) H_{1}(Y_{j}(\bar{x})) H_{1}(Y_{k}(\bar{x})) H_{4}(f_{\ell
}(y(\phi )))\right] \sin \phi d\phi\\
&\;\;+16\pi^{2} \lambda_{\ell}
\sum_{i \ne j \ne k \ne l} \{\lim_{\varepsilon \rightarrow 0} q_{ijkl}^{\varepsilon }(I)
\} \int_{0}^{\pi /2}\mathbb{E}\left[ H_{1}(Y_{i}(\bar{x})) H_{1}(Y_{j}(\bar{x})) H_{1}(Y_{k}(\bar{x})) H_{1}(Y_{l}(\bar{x})) H_{4}(f_{\ell
}(y(\phi )))\right] \sin \phi d\phi.
\end{align*}
We observe that the asymptotic behaviour of $\text{Cov}(\mathcal{N}_{\ell
}^{c}(I),h_{\ell ;4})$ is dominated by three terms corresponding to
\begin{equation*}
\int_{0}^{\pi /2}\mathbb{E}\left[ H_{4}(Y_{2}(\bar{x}))H_{4}(f_{\ell
}(y(\phi )))\right] \sin \phi d\phi, \;\;\; \int_{0}^{\pi /2}\mathbb{E}
\left[ H_{4}(Y_{5}(\bar{x}))H_{4}(f_{\ell }(y(\phi )))\right] \sin \phi d\phi,
\end{equation*}
and
\begin{equation*}
\int_{0}^{\pi /2}\mathbb{E}\left[ H_{2}(Y_{2}(\bar{x}))H_{2}(Y_{5}(\bar{x}
))H_{4}(f_{\ell }(y(\phi )))\right] \sin \phi d\phi.
\end{equation*}
The computation of these leading covariances is given in Lemma \ref{Lem5.1-5.2-5.3}. All the remaining terms in $\text{Cov}(\mathcal{N}_{\ell }^{c}(I),h_{\ell;4})$ are shown to be $O(\ell^{-2})$ or smaller in Lemma \ref{lem5.4,5.5,5.6}.
From Proposition \ref{coefficients} we know that
\begin{equation*}
k_{2}(I)=
\frac{3}{8\pi } \cal{I}_0(I) , \hspace{1cm} k_{5}(I)=
\frac{3}{8\pi }\cal{I}_0(I)-\frac{2^3}{2^5\pi} \cal{I}_2(I)+ \frac{2^6}{3^2 2^9 \pi} \cal{I}_4(I),
\end{equation*}
\begin{equation*}
 h_{25}(I)
 =\frac{1}{2^3 \pi} \cal{I}_0(I)-\frac{2^3}{2^6 3 \pi} \cal{I}_2(I).
\end{equation*}
Substituting and after some straightforward algebra, one obtains
\begin{align*}
\text{Cov}&(\mathcal{N}_{\ell }^{c}(I),h_{\ell ;4}) \\&=\lambda _{\ell }\left\{ 4 \pi
^{2}h_{25}(I)4!\frac{3}{\pi ^{2}}\frac{\log \ell }{\ell ^{2}}+\frac{2}{3}\pi
^{2}k_{2}(I)4!2^{2}\frac{3}{\pi ^{2}}\frac{\log \ell }{2\ell ^{2}}+\frac{2}{3}
\pi ^{2}k_{5}(I)4!3^{2}\frac{3}{\pi ^{2}}\frac{\log \ell }{2\ell ^{2}}+O(\ell
^{-2})\right\} \\
&=\lambda_{\ell} \frac{\log \ell }{\ell^2} 4! \frac{1}{\pi} \frac{1}{2^9} \bigg\{ 51\cdot2^6 \cal{I}_0(I)-11\cdot 2^7 \cal{I}_2(I)+2^6 \cal{I}_{4}(I) \bigg\} +O(1). 
\end{align*}
The statement of Proposition \ref{lem2} follows recalling the definition of $\cal{F}_\ell(I)$ and Corollary \ref{lem1}. 

\subsection{Proof of Proposition \ref{coefficients}: evaluation of the Projection Coefficients $h_{52}(I),k_{2}(I),k_{5}(I)$}\label{EvaCoeff}

In this section we evaluate the three projection coefficients in the
Wiener-chaos expansion which are required in Proposition \ref{lem2}. Let us recall first the following simple result: assuming $Y$ standard Gaussian 
	\begin{equation*}
	\varphi_{i}:=\lim_{\varepsilon \to 0} \mathbb{E}[H_{i}(Y)\delta_{\varepsilon}(\mu_1 Y)]=
	\begin{cases}
	\frac{1}{\sqrt{2\pi }\mu_1} & i=0, \\
	0 & i=1, \\
	-\frac{1}{\sqrt{2\pi }\mu_1} & i=2, \\
	\frac{3}{\sqrt{2\pi }\mu_1} & i=4,
	\end{cases}
	\end{equation*}
	indeed, for example, 
	\begin{equation*}
	\lim_{\varepsilon \to 0} \mathbb{E}[H_{4}(Y)\delta_{\varepsilon} (\mu_1 Y)]=\lim_{\varepsilon \to 0} \mathbb{E}[(Y^{4}-6Y^{2}+3) \delta_{\varepsilon} (\mu_1 Y)]=\frac{3}{\sqrt{2\pi }\mu_1},
	\end{equation*}
	since
	\begin{equation*}
	\lim_{\varepsilon \to 0} \mathbb{E}[Y^{n} \delta_{\varepsilon}(\mu_1 Y)]=\lim_{\varepsilon \to 0}\int_{-\infty }^{\infty }y^{n}\delta_{\varepsilon} (\mu_1 y)\frac{1}{\sqrt{2\pi }}e^{-\frac{y^{2}}{2}}dy=
	\begin{cases}
	\frac{1}{\sqrt{2\pi }\mu_1} & n=0, \\
	0 & n=1,2,3\dots
	\end{cases}
	\end{equation*}
	This allows us to write $k_{2}(I)$, $k_{5}(I)$ and $k_{52}(I)$ as follows: 
	\begin{align*}
	k_{2}(I)&=\lim_{\varepsilon \rightarrow 0 }k^{\varepsilon}_{2}(I) \\
	&=  \lambda_{\ell } \; \mathbb{E}\left[ \left\vert \frac{\mu_{3}\mu_{5}}{\lambda _{\ell }^{2}}Y_{3}Y_{5}+\frac{\mu_{2}\mu_{3}}{\lambda _{\ell }^{2}}Y_{3}^{2}-\frac{\mu_{4}^{2}}{\lambda
	_{\ell }^{2}}Y_{4}^{2} \right\vert \id_{\{ \frac{\mu_2+\mu_3}{\lambda_\ell} Y_3+\frac{\mu_5}{\lambda_\ell} Y_5 \in I \}} \right] \varphi _{0}\,\varphi _{4} \\
	&= \frac{3}{\pi }  \mathbb{E}\left[ \left\vert \frac{1}{2\sqrt{2}}Y_{3}Y_{5}+
	\frac{1}{8}Y_{3}^{2}-\frac{1}{8}Y_{4}^{2}\right\vert \id_{\{ \frac{\sqrt 2}{\sqrt 3}  Y_3+ \frac{1}{\sqrt 3}  Y_5 \in I \}} \right] +O(\ell^{-1}),
	\end{align*}
	\begin{align*}
	k_{5}(I)&=\lim_{\varepsilon \rightarrow 0 }k^{\varepsilon}_{5} (I)\\
	&= \lambda_{\ell } \;\mathbb{E}\left[ \left\vert \frac{\mu_{3}\mu_{5}}{\lambda _{\ell }^{2}}Y_{3}Y_{5}+\frac{\mu
	_{2}\mu_{3}}{\lambda _{\ell }^{2}}Y_{3}^{2}-\frac{\mu_{4}^{2}}{\lambda
	_{\ell }^{2}}Y_{4}^{2} \right\vert \id_{\{ \frac{\mu_2+\mu_3}{\lambda_\ell} Y_3+\frac{\mu_5}{\lambda_\ell} Y_5 \in I \}}  H_{4}(Y_{5})
	\right] \varphi _{0}^{2} \\ &= \frac{1}{\pi }\;\mathbb{E}\left[ \left\vert \frac{1}{2\sqrt{2}}
	Y_{3}Y_{5}+\frac{1}{8}Y_{3}^{2}-\frac{1}{8}Y_{4}^{2}\right\vert \id_{\{ \frac{\sqrt 2}{\sqrt 3}  Y_3+ \frac{1}{\sqrt 3}  Y_5 \in I \}}  H_{4}(Y_{5})
	\right]+O(\ell^{-1}),
	\end{align*}
	and
	\begin{align*}
	h_{52}(I)&=\lim_{\varepsilon \rightarrow 0}h^{\varepsilon}_{25} (I)\\
	&=  \lambda_{\ell }\; \mathbb{E}\left[ \left\vert \frac{\mu_{3}\mu_{5}}{\lambda _{\ell }^{2}}Y_{3}Y_{5}+\frac{\mu
	_{2}\mu_{3}}{\lambda _{\ell }^{2}}Y_{3}^{2}-\frac{\mu_{4}^{2}}{\lambda
	_{\ell }^{2}}Y_{4}^{2} \right\vert \id_{\{ \frac{\mu_2+\mu_3}{\lambda_\ell} Y_3+\frac{\mu_5}{\lambda_\ell} Y_5 \in I \}}  H_{2}(Y_{5})\right] \varphi
	_{0}\,\varphi _{2} \\
	&=  - \frac{1}{\pi } \mathbb{E}\left[ \left\vert \frac{1}{2\sqrt{2}}Y_{3}Y_{5}+
	\frac{1}{8}Y_{3}^{2}-\frac{1}{8}Y_{4}^{2}\right\vert \id_{\{ \frac{\sqrt 2}{\sqrt 3}  Y_3+ \frac{1}{\sqrt 3}  Y_5 \in I \}}  H_{2}(Y_{5})\right] +O(\ell^{-1}).
	\end{align*}
	Let us introduce the change of variables
	\begin{equation*}
	Z_{1}=\sqrt{3}Y_{3},\hspace{1cm}Z_{2}=Y_{4},\hspace{1cm}Z_{3}=\frac{\sqrt{8}
	}{\sqrt{3}}Y_{5}+\frac{1}{\sqrt{3}}Y_{3},
	\end{equation*}
	so that $(Z_1, Z_2, Z_2)$ is a centred Gaussian vector with covariance matrix $\tilde{c}_{\infty }$
	and we can write
	\begin{align*}
	\frac{1}{2\sqrt{2}}Y_{3}Y_{5}+\frac{1}{8}Y_{3}^{2}-\frac{1}{8}Y_{4}^{2}=\frac{1}{8} (Z_{1}Z_{3}- Z_{2}^{2}), \hspace{0.5cm} \frac{\sqrt 2}{\sqrt 3}  Y_3+ \frac{1}{\sqrt 3}  Y_5  = \frac{1}{\sqrt{8}}Z_1+\frac{1}{\sqrt{8}}Z_3. \end{align*}
	Hence 
	\begin{equation*}
	\begin{split}
	k_{2}(I) &=\frac{3}{\pi } \mathbb{E}\left[ \frac{1}{8}\left\vert Z_1Z_3-Z_2^2 \right\vert \id_{\{ \frac{1}{\sqrt{8}}(Z_1+Z_3) \in I \}}\right] \\&
	=\frac{3}{\pi } \mathbb{E}\left[ \mathbb{E}\left[ \frac{1}{8}\left\vert Z_1Z_3-Z_2^2 \right\vert \id_{\{ \frac{1}{\sqrt{8}}(Z_1+Z_3) \in I \}}\bigg\vert Z_1+Z_3=\sqrt{8}t\right] \right]
	\\&
	= \frac{3}{8\pi} \int_I p_0(t) \,dt=\frac{3}{8\pi}  \cal{I}_0(I),
	\end{split}
	\end{equation*}
as claimed. Similarly,
	\begin{align*}
	h_{52}(I)
	& =-\frac{1}{\pi }\frac{1}{8}\mathbb{E}\left[ \left\vert
	Z_{1}Z_{3}-Z_{2}^{2}\right\vert \id_{\{ \frac{1}{\sqrt{8}}(Z_1+Z_3) \in I \}} H_{2}\left( \frac{1}{\sqrt{8}\sqrt{3}}
	(3Z_{3}-Z_{1})\right) \right]\\
	& =-\frac{1}{\pi }\frac{1}{8}\mathbb{E}\left[ \left\vert
	Z_{1}Z_{3}-Z_{2}^{2}\right\vert \id_{\{ \frac{1}{\sqrt{8}}(Z_1+Z_3) \in I \}}\left( \frac{1}{\sqrt{8}\sqrt{3}}
	(3Z_{3}-Z_{1})\right) ^{2}\right] \\&+\frac{1}{\pi }\frac{1}{8}\mathbb{E}\left[
	\left\vert Z_{1}Z_{3}-Z_{2}^{2}\right\vert \id_{\{ \frac{1}{\sqrt{8}}(Z_1+Z_3) \in I \}}\right] \\&
	=-\frac{1}{\pi } \frac{1}{8^2 3} \int_I \sqrt{8} \mathbb{E}\left[ \left\vert Z_{1}Z_{3}-Z_{2}^{2}\right\vert (3Z_3-Z_1)^2 \bigg|Z_1+Z_3=\sqrt{8}t\right]\phi_{Z_1+Z_3}(\sqrt{8}t)\,dt \\&+\frac{1}{8 \pi } \int_I \mathbb{E}\left[
	\left\vert Z_{1}Z_{3}-Z_{2}^{2}\right\vert \bigg| Z_1+Z_3=\sqrt{8}t\right] \phi_{Z_1+Z_3}(\sqrt{8}t)\,dt
	\\&=-\frac{2^3}{3\cdot 2^{6} \pi}\mathcal{I}_{2}(I)+\frac{1
	}{2^{3 }\pi }\mathcal{I}_{0}(I),
	\end{align*}
	and
	\begin{align*}
	k_{5}(I)
	& =\frac{1}{\pi }\;\frac{1}{8}\mathbb{E}\left[ \left\vert
	Z_{1}Z_{3}-Z_{2}^{2}\right\vert  \id_{\{ \frac{1}{\sqrt{8}}(Z_1+Z_3) \in I \}}H_{4}\left( \frac{1}{\sqrt{8}\sqrt{3}}
	(3Z_{3}-Z_{1})\right) \right]\\&
	= \frac{1}{8^3 3^2\pi} \int_I \sqrt{8} \mathbb{E}[ |Z_1Z_3-Z_2^2| (3Z_3-Z_1)^4 \bigg| Z_1+Z_3=\sqrt{8}t] \phi_{Z_1+Z_3}(\sqrt{8}t)\,dt\\&
	-\frac{1}{8\pi 4} \int_I \sqrt{8}\mathbb{E}[ |Z_1Z_3-Z_2^2| (3Z_3-Z_1)^2 \bigg| Z_1+Z_3=\sqrt{8}t] \phi_{Z_1+Z_3}(\sqrt{8}t)\,dt
	\\&+\frac{3}{8\pi } \int_I \sqrt{8}\mathbb{E}[ |Z_1Z_3-Z_2^2| \bigg| Z_1+Z_3=\sqrt{8}t] \phi_{Z_1+Z_3}(\sqrt{8}t)\,dt \\& =\frac{2^6}{2^{9}\cdot 3^{2} \pi}\mathcal{I}_{4}(I)-
	\frac{2^3}{2^{5} \pi}\mathcal{I}_{2}(I)+\frac{3}{2^{3} \pi}\mathcal{I}_{0}(I).
	\end{align*}
	
\subsection{Proof of Proposition \ref{17:37}: terms with odd index Hermite polynomials} \label{oddt}

The terms in the $4$-th chaos formula \eqref{4chaos} with odd index Hermite polynomials produce in $\text{Cov}(\mathcal{N}_{\ell }^{c}(I),h_{\ell;4})$ terms of order $O(\ell^{-2})$ and terms equal to zero, in fact, recalling that for $a$ odd we have
	\begin{equation*}
	\lim_{\varepsilon \to 0} \mathbb{E}[H_{a}(Y)\delta_{\varepsilon}(\mu_1 Y)]=0,
	\end{equation*}
we immediately see that the coefficients $g_{ij}(I)$ with $i, j = 1, 2$ are all equal to zero. For the coefficients $g_{ij}(I)$ with $i=4$ or $j=4$, we observe that the expectation with respect to the random variable $Y_4$ vanishes since it is expressed as the integral of an odd function. The proof of the last two points of the statement is similar.

\section{Proof of Proposition \ref{EPC4}} \label{EPCsection}
We define the approximating sequence 
$$\mathcal{L}_{0,\varepsilon}(u,\ell)=\int_{\mathbb{S}^2} (\det \nabla^2 f_{\ell}(x) ) \id_{\left\{ f_{\ell}(x) \ge u\right\}} \delta(\nabla f_{\ell}(x)) d x.$$
Under the assumptions of \cite[Lemma 4]{CM2018}, we have that 
$$\mathcal{L}_0(u,\ell) =\lim_{\varepsilon \to 0} \mathcal{L}_{0,\varepsilon}(u,\ell)$$
where the convergence holds both $\omega-a.s.$ and in $L^2(\Omega)$. The proof of Proposition \ref{EPC4} follows the same lines of the proof of Proposition \ref{lem2} with the only difference that the relevant projection coefficients are now 
\begin{align*}
& h_{25}(u)\\
&= \lim_{\varepsilon \to 0}\mathbb{E}\left[ \left( 
\frac{\mu _{3}\mu_{5}}{\lambda_\ell ^{2}}Y_{3}Y_{5}+\frac{\mu
	_{2}\mu_{3}}{\lambda_\ell ^{2}}Y_{3}^{2}-\frac{\mu _{4}^{2}}{\lambda_\ell ^{2}
}Y_{4}^{2}\right)  \id_{\left\{ \frac{\mu_2+\mu_3}{\lambda_\ell} Y_3+ \frac{\mu_5}{\lambda_\ell} Y_5 \le u \right\}} \delta_{\varepsilon} (Y_{1},Y_{2})H_{2}(Y_{2})H_{2}(Y_{5})\right], 
\end{align*}
and, for $i=2,5$, 
\begin{align*}
& k_{i}(u )=  \lim_{\varepsilon \to 0} \mathbb{E}\left[ \left( 
\frac{\mu_{3}\mu _{5}}{\lambda_\ell ^{2}}Y_{3}Y_{5}+\frac{\mu_{2}\mu _{3}}{\lambda_\ell ^{2}}Y_{3}^{2}-\frac{\mu_{4}^{2}}{\lambda_\ell ^{2}
}Y_{4}^{2}\right)  \id_{\left\{ \frac{\mu_2+\mu_3}{\lambda_\ell} Y_3+ \frac{\mu_5}{\lambda_\ell} Y_5 \le u \right\}} \delta_{\varepsilon}(Y_{1},Y_{2})H_{4}(Y_{i})\right]. 
\end{align*}
After a long series of calculations that we do not include here for brevity sake,  we obtain the following explicit expression for the projection coefficients: 
\begin{align*}
h_{25}(u)&=-\frac{1}{24 \pi} u (u^2+1) \frac{e^{-u^2/2}}{\sqrt{2 \pi}},\\
k_{2}(u)&=\frac{3}{8 \pi} u \frac{e^{-u^2/2}}{\sqrt{2 \pi}},\\
 k_{5}(u)&=\frac{1}{27 \pi} u (-9-2u^2+u^4) \frac{e^{-u^2/2}}{\sqrt{2 \pi}}. 
\end{align*}
The statement immediately follow in the nodal case $u=0$.

\end{document}